\documentclass[11pt]{article}

\usepackage[a4paper,margin=1in]{geometry}
\usepackage{amsmath,amssymb,amsfonts,amsthm}
\usepackage[showonlyrefs]{mathtools}
\usepackage{mathrsfs}
\usepackage{graphicx}
\usepackage{booktabs}
\usepackage{multirow}
\usepackage{algorithm}
\usepackage{algorithmicx}
\usepackage{algpseudocode}
\usepackage{enumitem}
\usepackage[title]{appendix}
\usepackage{xcolor}
\usepackage{textcomp}
\usepackage{manyfoot}
\usepackage{listings}
\usepackage[colorlinks=true, allcolors=blue]{hyperref}
\usepackage[capitalize]{cleveref}
\usepackage[
  style=numeric,      
  giveninits=true,
  sorting=none,
  natbib=true,        
  doi=false,
  url=false,
  isbn=false,
  eprint=false
]{biblatex}
\DeclareFieldFormat[article,incollection,inproceedings,misc]{title}{#1}
\renewbibmacro{in:}{}
\renewbibmacro*{volume+number+eid}{\printfield{volume}\printfield[parens]{number}}
\DeclareFieldFormat{pages}{#1}

\DeclareFieldFormat[article,incollection,inproceedings,unpublished]{title}{#1}
\renewbibmacro{in:}{}
\renewbibmacro*{volume+number+eid}{\printfield{volume}\printfield[parens]{number}}
\renewbibmacro*{journal+issuetitle}{\printfield{journaltitle}\setunit*{\addcomma\space}}
\DeclareFieldFormat[article,inproceedings,incollection]{pages}{#1}
\DeclareBibliographyDriver{article}{%
  \usebibmacro{bibindex}\usebibmacro{begentry}%
  \usebibmacro{author}%
  \setunit{\labelnamepunct}\newblock
  \usebibmacro{title}\newunit\newblock
  \usebibmacro{journal+issuetitle}%
  \usebibmacro{volume+number+eid}%
  \addcolon\space\printfield{pages}%
  \addcomma\space\printfield{year}%
  \usebibmacro{finentry}}
\DeclareBibliographyDriver{incollection}{%
  \usebibmacro{bibindex}\usebibmacro{begentry}%
  \usebibmacro{author}\setunit{\labelnamepunct}\newblock
  \usebibmacro{title}\newunit\newblock
  \printfield{booktitle}\setunit{\addcomma\space}%
  \printlist{publisher}\setunit{\addcomma\space}%
  \printfield{pages}\addcomma\space\printfield{year}%
  \usebibmacro{finentry}}
\DeclareBibliographyAlias{inproceedings}{incollection}
\DeclareBibliographyDriver{unpublished}{%
  \usebibmacro{bibindex}\usebibmacro{begentry}%
  \usebibmacro{author}\setunit{\labelnamepunct}\newblock
  \usebibmacro{title}\newunit\newblock
  \printfield{note}
  \setunit{\addcomma\space}
  \printfield{year}%
  \usebibmacro{finentry}}

\newcommand{\R}{\operatorname{\mathbb R}}
\newcommand{\eR}{\operatorname{\overline{\mathbb R}}}
\newcommand{\Fix}{\operatorname{\mathcal{F}}}
\newcommand{\rank}{\operatorname{rank}}
\newcommand{\Null}{\operatorname{Null}}
\newcommand{\range}{\operatorname{range}}
\newcommand{\prox}{\operatorname{prox}}
\newcommand{\refop}{\operatorname{ref}}

\newcommand{\cJ}{\operatorname{\mathcal J}}
\newcommand{\dist}{\operatorname{dist}}

\newcommand{\tr}{^\intercal}
\newcommand{\argmin}{\operatorname{argmin}}
\newcommand{\N}{\mathbb{N}}
\newcommand{\st}{\operatorname{st.}}
\newcommand{\kQ}{\kappa^+(Q)}
\newcommand{\barx}{\overline{x}}
\newcommand{\dmin}{\displaystyle \min}
\newcommand{\dom}{\operatorname{dom}}

\newtheorem{theorem}{Theorem}
\newtheorem{proposition}[theorem]{Proposition}
\newtheorem{lemma}{Lemma}
\newtheorem{corollary}{Corollary}
\newtheorem{definition}{Definition}
\newtheorem{example}{Example}
\newtheorem{remark}{Remark}

\title{Linear Convergence and Error Bounds for Optimization Without Strong Convexity}

\author{Kira van Treek\thanks{Econometrics and Operations Research, Tilburg University,
Warandelaan 2, 5037 AB Tilburg, The Netherlands. 
Email: \texttt{k.vantreek@tilburguniversity.edu}}
\and Javier F. Pe{\~n}a\thanks{Tepper School of Business, Carnegie Mellon University,
5000 Forbes Ave, Pittsburgh, PA 15213, United States.
Email: \texttt{jfp@andrew.cmu.edu}}
\and Juan C. Vera\thanks{Econometrics and Operations Research, Tilburg University,
Warandelaan 2, 5037 AB Tilburg, The Netherlands.
Email: \texttt{j.c.veralizcano@tilburguniversity.edu}}
\and Luis F. Zuluaga\thanks{Industrial and Systems Engineering, Lehigh University,
27 Memorial Drive West, Bethlehem, PA 18015, United States.
Email: \texttt{luis.zuluaga@lehigh.edu}}}

\date{} 

\begin{document}
\maketitle
\begin{abstract}
Many optimization algorithms—including gradient descent, proximal methods, and operator splitting techniques—can be formulated as fixed-point iterations (FPI) of continuous operators. When these operators are averaged, convergence to a fixed point is guaranteed when one exists, but the convergence is generally sublinear. Recent results establish linear convergence of FPI for averaged operators under certain conditions. However, such conditions do not apply to common classes of operators, such as those arising in piecewise linear and quadratic optimization problems. In this work, we prove that a local error-bound condition is both necessary and sufficient for the linear convergence of FPI applied to averaged operators. We provide explicit bounds on the convergence rate and show how these relate to the constants in the error-bound condition. Our main result demonstrates that piecewise linear operators satisfy local error bounds, ensuring linear convergence of the associated optimization algorithms. This leads to a general and practical framework for analyzing convergence behavior in algorithms such as ADMM and Douglas–Rachford in the absence of strong convexity. In particular, we obtain convergence rates that are independent of problem data for linear optimization, and depend only on the condition number of the objective for quadratic optimization.
\end{abstract}



\section{Introduction.}
A wide range of iterative numerical methods to solve optimization problems such as the gradient descent, the proximal point, proximal gradient,  the Douglas-Rachford and the alternating direction method of multipliers (ADMM), among many others, have in common that they can be stated by
iteratively applying an appropriate continuous operator $F:\R^n \to \R^n$ to a given starting point. Namely, starting with an initial guess $x_0 \in \R^n$,  these methods iteratively apply the operator $F$, to obtain the sequence $x_{k+1} = F(x_k)$, which is known as \textit{fixed-point iteration} (FPI). The aim of these methods is to obtain a sequence $x_k$, $k=0,1,\dots$ that converges to a point $x$ in the set of optimal solutions to the optimization problem of interest.
Notice that the continuity of $F$ implies that when such $x$ exists, it satisfies the \textit{fixed-point} equation
\[
 F(x) = x.
\]
We will denote the set of fixed-points associated to an operator $F$ as
\(
\Fix_{F}:= \{x \in \R^n: F(x) = x\}.
\)
When solving optimization problems, $F$ is chosen such that the set $\Fix_{F}$ is related to the set of optimal points for the optimization problem, in such fashion that any $x \in \Fix_F$ can be used to construct an optimal point.

\begin{example}[Gradient descent for unconstrained convex optimization] \label{ex:gradient}
Arguably one of the simplest examples of the use of FPI for optimization is the Gradient descent method.
Given a convex function $f: \R^n \to \R$, consider the unconstrained optimization problem
\begin{equation}
\label{eq:minf}
\min_{x \in \R^n} f(x).
\end{equation}
 A standard method to solve~\eqref{eq:minf} when $f$ is differentiable is the gradient descent algorithm, which consists of applying FPI to $F_{GD}(x) = x - \lambda \nabla f(x)$, where $\lambda$ is a small positive constant. Notice that by convexity, the set of fixed-points $\Fix_{F_{GD}}= \{x \in \R^n:F_{GD}(x) = x\} = \{x \in \R^n:\nabla f(x) = 0\}=\argmin_x f(x)$, the set of optimal points of~\eqref{eq:minf}.
\end{example}

Analogous to Example~\ref{ex:gradient}, in Table~\ref{tab:FPIs}, the fixed-point operator $F$ associated to the FPI formulation of optimization algorithms such as the proximal point, gradient projection, proximal gradient, Douglas-Rachford, Peaceman-Rachford, and ADMM  algorithms is explicitly given in terms of gradient, proximal, reflector, and projection operators. For details see \cref{sec:prelim}.

\begin{table}[htb]
\centering
{\small
\begin{tabular}{lll}
  \textbf{Algorithm}&   \textbf{Problem} &  \textbf{Fixed-Point Operator} \\
\toprule
Gradient descent &     \multirow{2}*{$\displaystyle \min_{x \in \R^n} f(x)$} & $F_{GD}(x):=x-\lambda \nabla f(x)$ \\
 Proximal point &  & $F_{P}(x):=\prox_{\gamma f}(x)$ \\
 \midrule
 Gradient projection &  $\displaystyle \min_{x \in S} f(x)$ & $F_{GP}(x):=\Pi_{S}(x-\lambda \nabla f(x))$ \\ 
 \midrule
 Proximal gradient &   \multirow{3}*{$\displaystyle \min_{x \in \R^n} f(x)+g(x)$}  & $F_{PG}(x):=\prox_{\gamma g}(x-\lambda \nabla f(x))$ \\
 Peaceman-Rachford & &$F_{PR}(x):= \refop_{\gamma g}(\refop_{\gamma f}(x))$\\
 Douglas-Rachford & &$F_{DR}(x):= (1 - \alpha) x  + \alpha\refop_{\gamma g}(\refop_{\gamma f}(x))$\\
 \midrule
 ADMM &
$
\begin{array}{cl}
\displaystyle \min_{x, y \in \R^n} & f(x)+g(y)\\
\st & Ax + By = c
\end{array}
$
&
$
\begin{array}{l}
 F_{ADMM}(x) :=  \tfrac{1}{2}x  + \tfrac{1}{2}\refop_{\gamma \tilde g}(\refop_{\gamma \tilde f}(x))\\ 
  \hfill \text{with } \tilde{g}(x) := \min_{z \in \R^n} \{ g(z): Bz = x\},\\
  \hfill \tilde{f}(x) := \min_{z \in \R^n} \{ f(z): c-Az = x\} \\
 \end{array}
 $\\
\bottomrule
\end{tabular}
}
\caption{Fixed-point operators associated to optimization algorithms in terms of gradient, proximal, reflector, and projection operators, where $f, g \in \R^n \to \R\cup \{+\infty\}$ with $f$ being differentiable in the cases where $\nabla f$,  the gradient of $f$,  is used, $\gamma > 0$, $S \subseteq \R^n$, $\alpha \in (0,1)$, $A, B \in \R^{m \times n}$, $c \in \R^m$.}
\label{tab:FPIs}
\end{table}

Also, analogous to \cref{ex:gradient}, the optimal solutions of all the algorithm-problem pairs in \cref{tab:FPIs} can be written in terms of the set of fixed-points of their associated fixed-point operator. This fact is formally stated and used later in \cref{sec:FPvsOP} (\cref{lem:fix.p.opt.con}).

Since the convergence of FPI is asymptotic, the iterative process is halted after a finite number of iterations \( k \), once an iterate that is sufficiently close to a fixed point is achieved. This is determined when \( \dist(x_k, \Fix_F) \), the (Euclidean) distance  to a fixed-point,  is sufficiently small. Since this distance cannot be directly computed, a practical stopping criterion is to ensure that the residual, given by \( \|F(x_k) - x_k\| = \|x_{k+1} - x_k\| \), is sufficiently small. Thus, it is crucial to demonstrate that the operator \( F \) satisfies a \emph{(Lipschitz) error-bound} condition that relates the distance to the fixed-points to the residual. This relationship ensures that the proposed stopping criterion for the iterative process is sound. If the operator does not satisfy such an error bound, the algorithm cannot guarantee a near-optimal solution after stopping.

Contracting operators satisfy a global error-bound condition and the corresponding algorithm converges linearly. Under strong convexity and smoothness assumptions all the operators in \cref{tab:FPIs} are contracting (see, \citep{giselsson2016linear} and \cref{sec:prelim}). Indeed, the proofs of the linear convergence of these algorithms under such assumptions, which are scattered around in the literature~\citep[see,][and references therein]{davis2017faster, giselsson2016linear, moursi2019douglas, abbaszadehpeivasti2022exact, karimi2016linear}, implicitly (and sometimes explicitly) use the contracting property of the corresponding operator. 

 It may seem logical to concentrate on designing contracting operators for optimization. However, an important consideration arises: Banach's fixed-point theorem ensures the existence of a unique fixed point. Consequently, when multiple optimizers are present, it is reasonable to expect that the corresponding operator will not be contractive, as we should anticipate having more than one fixed point in such cases. In particular, in the nonstrongly convex case, it is important to consider operators beyond the contracting ones.

A weakened version of contraction is nonexpansion. An operator $N$ is said to be {\em nonexpansive} if for all $x,y\in\R^n$,
\begin{equation}\label{def:nonexpansive}
 \|N(x) - N(y)\| \le \|x-y\|.
\end{equation}
Nonexpansive operators might not have any fixed points (e.g. translations); and when they do, the FPI sequence might not be convergent (e.g. rotations). A standard technique for finding fixed points of nonexpansive operators is to use the Krasnosel'ski\u{\i}–Mann iteration, which involves applying the FPI to an {\em averaging} of the nonexpansive operator.

Given $0 < \alpha < 1$, an operator $F$ is said to be {\em $\alpha$-averaged}, or simply {\em averaged}, if it is defined as: \[ F(x) = (1 - \alpha)x + \alpha N(x) \] for some nonexpansive operator $N$. It is important to note that the set of fixed points remains unchanged by the averaging operation; that is, $\Fix_{F} = \Fix_{N}$. 
As stated in \cref{thm:avgcon} below, the FPI of averaged operators always converges to a fixed point if one exists. However, in general, this convergence is only sublinear, and error bounds are not guaranteed.

\begin{example}[Gradient descent is an averaged operator]
\label{ex:gradient2} Consider the gradient descent operator $F_{GD}(x) = x -\lambda \nabla f(x)$. By definition, $F_{GD}$ is non-expansive if and only if for all $x$ and $y$,
\[
\|x-y\|^2 \ge \|F_{GD}(x)-F_{GD}(y)\|^2 =  \|x-y\|^2  - 2\lambda(x-y)\tr(\nabla f(x) - \nabla f(y)) + \lambda^2\|\nabla f(x) - \nabla f(y)\|^2.
\]
Which is equivalent to,
\[
(x-y)\tr(\nabla f(x) - \nabla f(y)) \ge \frac{\lambda}{2}\|\nabla f(x) - \nabla f(y)\|^2,
\]
which is in turn equivalent to $f$ being $\tfrac{\lambda}{2}$-smooth. 
Now, we turn our attention to when is $F_{GD}$ an averaged operator. Given $\alpha \in (0,1)$.  We write the gradient descent operator as $F_{GD}(x) =x -\lambda \nabla f(x) = (1-\alpha) x + \alpha F_{GD}^{\alpha}(x)$, where $F_{GD}^{\alpha}(x) := (x - \frac{\lambda}{\alpha} \nabla f(x))$. Then $F_{GD}$ is $\alpha$-averaged if and only if $F_{GD}^{\alpha}(x)$ is non-expansive, which is equivalent to $f$ being $\tfrac{\lambda}{2\alpha}$-smooth.
\end{example}

\begin{theorem}[Convergence of averaged operators]\label{thm:avgcon}
Let $0 < \alpha < 1$ and $F:\R^n\to\R^n$ be an $\alpha$-averaged operator, such that $F$ has fixed-points. Given $x_0\in\R^n$, let $x_{k+1}=F(x_k)$ for $k=0,1,\dots$. Then
\begin{enumerate}[label=(\roman*)]
    \item $x^*= \lim_{k \to \infty} x_k$ exists and $F(x^*) = x^*$~\citep[][Thm.~5.14]{baus2017convex}.
    \item The rate of convergence of $\|x_{k+1}-x_k\|\to 0$ is of order $O(1/\sqrt{k})$ \citep{baillon1996the}.\end{enumerate}
\end{theorem}

In practice, parameters are chosen so that the fixed-point operators associated with various optimization algorithms are averaged. As discussed in detail in Section~\ref{sec:prelim}, all the fixed-point operators associated with the algorithms listed in Table~\ref{tab:FPIs} are averaged, except for the Peaceman-Rachford operator. As elaborated in Section~\ref{sec:PR}, the Douglas-Rachford operator is the averaged version of the Peaceman-Rachford operator.

\cref{thm:avgcon} is a fundamental result applicable to a wide array of optimization algorithms. However, it merely offers sublinear convergence and does not provide any specific error bounds. Recently, this result has been strengthened by showing that FPI of any averaged operator converges linearly to its fixed-points when the operator satisfies an error-bound condition, paving the way for generic techniques covering large classes of algorithms and/or models. In~\citep{banjac2018tight} this type of result is established under a
{\em global} error-bound condition. However, this global error bound condition limits the applicability of the result, as piecewise linear operators typically satisfy only local error bound conditions. Such operators commonly arise in important classes of optimization problems that are not strongly convex, including linear and quadratic optimization. In~\citep{aspelmeier2016local}, the global error-bound condition is weakened to a local error-bound condition. However, the result requires the existence of a unique fixed-point which again limits the applicability of the result to the cases discussed above.

\subsection*{Contributions.}
Here, we begin by enhancing the results in~\citep{aspelmeier2016local, banjac2018tight} by using a local error-bound condition similar to the one used in~\citep{aspelmeier2016local} with no extra assumptions beyond the averageness of the operator. We provide both upper and lower bounds for the associated linear rate of convergence, as well as bounds not only for the convergence in {\em distance} form (which coincides with the one given in~\citep[][Thm. 1]{banjac2018tight}), but also in {\em sequence} form. More specifically, we show  (\cref{thm.FPlinconvergence} and \cref{prop.lcEB})
that for an averaged operator $F$, the linear convergence of the associated FPI is equivalent to the following ({\em local}) error-bound condition for a suitable $R > 0$:
\begin{equation} \label{eq.error.bound.F}
\begin{gathered}
 \text{There exists $K_F>0$ such that }\\
 \dist(x,\Fix_F) \le R \text{~implies~} \dist(x,\Fix_F) \le K_F\cdot\|F(x) - x\|.
 \end{gathered}
\end{equation}
In particular, if the error-bound condition~\eqref{eq.error.bound.F} holds for some $\alpha$-averaged operator $F$ with
$\mathcal{F}_F \neq \emptyset$, then the sequence of points $x_{k+1}=F(x_k)$ for $k=0,1,\dots$,  generated by the FPI of $F$ starting from some initial iterate $x_0 \in \R^n$, converges linearly to $\mathcal{F}_F$ at a rate that can be bound in terms of the constant $K_F$ in the error-bound condition~\eqref{eq.error.bound.F} as follows
\[\dist(x_{k+1},\Fix_F) \le \rho\dist(x_{k},\Fix_F), \text{ and }  \|x_{k+1}-\barx\| \le \left(1 - \tfrac{1+\rho}{2}(1-\rho)^2\right)^{\frac{1}{2}}  \|x_k-\barx\|,\]
where $\rho = \left (1-\tfrac{1-\alpha}{\alpha K_F^2}\right )^\frac{1}{2}$. Above, in the left, the linear convergence is stated in {\em distance} form, while in the right, it is given in {\em sequence} form. The equivalence between the error-bound condition and the linear convergence of averaged operators can be exploited (\cref{prop.lcEB}) to obtain the following upper and lower bound relationships between the constants defining the error-bound condition and the linear convergence of FPI for an $\alpha$-averaged operator $F$:
\[1-\tfrac 1 {\tilde{K}_F}\leq \tilde{\rho}_F \leq \sqrt{1- \tfrac{1-\alpha}{\alpha \tilde{K}_F^2}} \text{ and }\sqrt{\tfrac{1-\alpha}{\alpha (1-\tilde{\rho}_F^2)}}\leq \tilde{K}_F\leq \tfrac 1{1-\tilde{\rho}_F},\]
where ${\tilde{\rho}_F}$ and $\tilde{K}_F$ denote the tightest constants such that the linear convergence and the error-bound condition respectively hold.

With these results established, our main result (\cref{thm:HpieceWiseLin}) is to show that a wide range of optimization algorithms, when applied to a broad class of optimization problems, satisfy the error-bound condition~\eqref{eq.error.bound.F} and therefore, achieve linear convergence. This result provides a powerful and generic framework for proving the linear convergence of optimization algorithms without relying on strong convexity assumptions. The algorithms we consider encompass, but are not limited to, gradient descent, proximal point methods, proximal gradient, Douglas-Rachford, and ADMM. Our analysis includes various fundamental types of optimization problems, notably linear optimization (LO) and convex quadratic optimization (QO). Furthermore, we cover least squares and its adaptations, such as LASSO and elastic net, along with support vector machines (SVM).
 
The error-bound condition guarantee presented in \cref{thm:HpieceWiseLin} offers a generic framework for estimating the error-bound condition parameter, which in turn determines the associated linear convergence rate. This is achieved by computing bounds on an appropriate {\em relative} Hoffman constant~\citep{pena2021new}, which, as its name suggests, is closely related to the well-established Hoffman constant~\citep{hoffman1952on}. 

To highlight the significance of this approach, we refine our findings for the ADMM and Douglas-Rachford algorithms to establish error bounds and linear convergence rates in their application to LO (\cref{sec:LP}) and QO (\cref{sec:QP}) problems. Both the Douglas-Rachford and ADMM algorithms are among the most extensively studied optimization methods. A considerable body of research on these algorithms focuses on instances where the function(s) being optimized exhibit strong convexity and/or smoothness~\citep[see, e.g.,][and references therein]{davis2017faster, giselsson2016linear, zamani2024the}. Other studies assume that an error-bound condition is satisfied~\citep[see, e.g.,][]{liang2016convergence, hong2017linear, yuan2020discerning, liu2018partial}.
  
For LO, when applying ADMM, we obtain an error-bound constant \( K_{F_{ADMM}} = 1 \), which translates into a linear rate of convergence \( \rho_{F_{ADMM}} \leq \frac{1}{2} \). This rate is independent of the problem's data (see Theorem \ref{thm:LP}). In contrast, related studies provide bounds on the convergence rate and error-bound constant for LO, which depend on the specific problem data. 
For instance, in \cite{eckstein1990alternating}, a particular splitting method for LO is proposed, which results in a linear convergence of ADMM that is influenced by the problem's data, specifically depending on the size of the problem. Similarly, in \cite{wang2017a}, by considering another splitting method for LO, a linear convergence of ADMM is obtained depending on the problem's data and the maximum radius of all iterates generated by the algorithm. The linear rates of convergence for LO discussed in \cite{mai2021fast} and \cite{liang2017local} can be characterized similarly. 

For Quadratic Optimization (QO), the derived error-bound constant and linear rate of convergence depend solely on the condition number \( \kQ = \lambda_{\max}(Q)/\lambda_{\min}^+(Q) \) of the data matrix~\( Q \) defining the problem's quadratic objective (see Theorem \ref{thm:QP}), where $\lambda_{\min}^+(\cdot)$ is the smallest positive eigenvalue. Specifically, we have \( K_{F_{ADMM}} \leq 6 \kQ \) and \( \rho_{F_{ADMM}} \leq 1 - \frac{1}{72 \kQ^2} \). These results stand in contrast to other studies where the convergence rates and error-bound constants for QO are either not explicitly provided or require a more complex analysis of the problem data. 
For example, in \cite{boley2013local}, linear convergence is demonstrated for ADMM; however, the actual rates are not stated, and the analysis necessitates that the problem has a unique optimal solution. In \cite{ghadimi2014optimal}, a linear convergence rate of ADMM for QO under a strong convexity assumption is established, with the rate depending on the spectral information of a matrix derived from the complete original problem data. Moreover, \cite{han2013local} presents a linear convergence rate of ADMM for QO without needing strong convexity, but expresses the rate in terms of the spectral properties of matrices that depend intricately on all the original problem data. Lastly,~\cite{patrinos2014douglas} obtains a sublinear rate of convergence for QO and establishes a linear convergence rate under strong convexity. Under the strong convexity assumptions, this last rate is improved in \cite{giselsson2016linear}.

We conclude this section by outlining the organization of the rest of the article. In Section~\ref{sec:equivalence}, we present the equivalence between the existence of error bounds for averaged operators and the linear convergence of their associated FPI. In Section~\ref{sec:EB4PWLopers}, we present our main result, namely,  that piecewise linear operators satisfy the error-bound condition.  Consequently, the corresponding FPI for these operators converges linearly. In Section~\ref{sec:prelim}, we present the fixed-point operators associated with several widely used optimization algorithms. In Section~\ref{sec:algs}, we show that for a broad class of optimization problems, these operators are indeed piecewise linear, ensuring that our earlier results are directly applicable. 
Finally, in Sections~\ref{sec:LP} and~\ref{sec:QP}, we will showcase our findings in the context of ADMM and DR, particularly for solving  LO and QO problems. We conclude in Section~\ref{sec:final} with some reflections on the implications of our work.

\section{Linear convergence and error bound equivalence.}
\label{sec:equivalence}
 In this section, we show that the FPI associated with averaged operators converges linearly if and only if it satisfies the error-bound condition~\eqref{eq.error.bound.F}. As mentioned in the introduction, our results and methods enhance similar work in~\cite{aspelmeier2016local, banjac2018tight}. 
As illustrated in Section~\ref{sec:opt}, these enhancements allow us to apply this equivalence result
to analyze the convergence of optimization algorithms for which the corresponding operator is piecewise linear without requiring strong convexity assumptions.

We begin by providing a useful result about averaged operators.
\begin{lemma}\label{cor:avgin1}
Let $\alpha\in (0,1)$ and let $F:\R^n\to\R^n$ be $\alpha$-averaged. Let $\hat{x}\in \Fix_F$ be a fixed-point of~$F$. Then
\begin{equation*}
\tfrac{1-\alpha}{\alpha}\|F(x)-x\|^2 +  \|F(x)-\hat{x}\|^2 \leq \|x-\hat{x}\|^2.
\end{equation*}
\end{lemma}
\proof{}
From the characterization of $\alpha$-averaged operators given in~\citep[][Prop. 4.35($iii$)]{baus2017convex}, it follows that
for all $x,\hat x\in\R^n$,
\[
\tfrac{1-\alpha}{\alpha}\|(x-F(x))-(\hat x-F(\hat x))\|^2 + \|F(x)-F(\hat x)\|^2 \le \|x-\hat x\|^2.\]
The result then follows by using the fact that $F(\hat{x}) = \hat{x}$.
\endproof

Next, we prove that for averaged operators, satisfying the error-bound condition~\eqref{eq.error.bound.F} implies linear convergence of the associated FPI.

\begin{theorem}[Error bound implies linear convergence] \label{thm.FPlinconvergence} Let $\alpha\in (0,1)$ and $F:\R^n\to\R^n$ be a continuous $\alpha$-averaged operator such that  $\Fix_F \ne \emptyset$ and the error-bound condition~\eqref{eq.error.bound.F} holds for some $R > 0$.
Let $x_0$ be given and let $x_{k+1} = F(x_k)$ for $k=0,1,\dots$.
Then $x_k \to \barx \in \Fix_F$ linearly. In particular, for any $k\ge0$ such that $\dist(x_k, \Fix_F) \le R$,
\begin{equation}
\label{eq:thm1-rho}
\dist(x_{k+1},\Fix_F) \le \rho \dist(x_{k},\Fix_F)
\text{ where }
\rho :=  \left (1-\tfrac{1-\alpha}{\alpha K_F^2}\right )^\frac{1}{2} \le 1-\tfrac{1-\alpha}{2\alpha K_F^2},
\end{equation}
and,
\begin{equation}
\label{eq:thm1-hatrho}
\|x_{k+1}-\barx\| \le \hat{\rho} \|x_k-\barx\|, \text{ where } \hat \rho :=  \left(1 - \tfrac{1+\rho}{2}(1-\rho)^2\right)^{\frac{1}{2}} \le  1-\tfrac{(1-\alpha)^2}{16\alpha^2 K_F^4}.
\end{equation}
\end{theorem}
\proof{}
By \cref{thm:avgcon}, the fixed-point iteration of averaged operators converge. Thus, we know the sequence $x_k$ converges to a point $\barx \in \Fix_F$.
Let $R > 0$ be such that the error-bound condition~\eqref{eq.error.bound.F} holds.
Let $k_0$ be the smallest $k \ge 0$ such that $\dist(x_{k},\Fix_F) \le  R$.
First, notice that for all $k \ge  k_0$,  we have  $\dist(x_{k},\Fix_F)\le R$ as averaged implies nonexpansive. 
Now, for each $k \ge k_0$, let
\begin{equation}
\label{eq:closest}
\barx_k \in \argmin_{x\in\Fix_F}\|x_k-x\|.
\end{equation}
Notice that $\barx_k \in \Fix_F$ is well defined as $\Fix_F$ is closed. $\barx_k$ satisfies $\|x_k - \barx_k\| = \dist(x_k,\Fix_F).$
Now, fix $k \ge k_0$. By \cref{cor:avgin1},
\begin{equation}\label{eq.averaged}
    \tfrac{1-\alpha}{\alpha}\|x_{k+1} -x_k\|^2 + \|x_{k+1}-\barx_k\|^2 \leq \|x_k-\barx_k\|^2.
\end{equation}
Moreover, by \eqref{eq.error.bound.F}
\begin{equation}\label{eq.error.bound.wk}
    \|x_k-\barx_k\|=\textup{dist}(x_k,\Fix_F)\le K_F \|x_{k+1} -x_k\|.
\end{equation}
It follows then that
\begin{align}
\nonumber
    \|x_{k+1}-\barx_{k+1}\|^2&\le \|x_{k+1}-\barx_k\|^2\\
    \nonumber
   & \le  \|x_k-\barx_k\|^2-\tfrac{1-\alpha}{\alpha}\|x_{k+1}-x_k\|^2 && (\text{by~\eqref{eq.averaged}})\\
   \label{eq:positive}
   &  \le  \left(1-\tfrac{1-\alpha}{\alpha K_F^2}\right)\|x_k-\barx_k\|^2. && (\text{by~\eqref{eq.error.bound.wk}})
\end{align}
Note that from~\eqref{eq:positive}, indeed $1-\smash{\tfrac{1-\alpha}{\alpha K_F^2}} \ge 0$ follows. Thus $\rho$ in \eqref{eq:thm1-rho} is well defined and $0 \le \rho < 1$. After taking square roots in~\eqref{eq:positive}, we obtain~\eqref{eq:thm1-rho}.

To prove~\eqref{eq:thm1-hatrho}, we use that for any $k = 0,1,\dots$ and $j=0,1,\dots$ we have 
\[
\|x_{k+j+1}-\barx_{k+j}\| = \|F(x_{k+j})-F(\barx_{k+j})\| \le \|x_{k+j}-\barx_{k+j}\|.
\]
Thus,
\begin{align*}
    \|\barx_{k+j+1}-\barx_{k+j}\|
    & \le  \|\barx_{k+j+1}-x_{k+j+1}\|+\|x_{k+j+1}-\barx_{k+j}\|\\
    & \le   \|\barx_{k+j+1}-x_{k+j+1}\| + \|x_{k+j}-\barx_{k+j}\| \\
    &  \le (\rho+1)\rho^{j}\|x_k-\barx_k\|. && (\text{after repeatedly using~\eqref{eq:thm1-rho}})
\end{align*}
In addition, $\barx_k\to\barx$, since $x_k\to\barx$ and $\|x_k-\barx_k\|\to0$. Thus,
\begin{align}\label{eq.wk.dist.bound}
\|x_{k}-\barx\|&\leq \|x_{k}-\barx_{k}\|+\sum_{j=0}^{\infty}\|\barx_{k+j+1}-\barx_{k+j}\|\nonumber\\
&\le \left(1 +(1+\rho)\sum_{j=0}^{\infty}\rho^{j}\right)\|x_k-\barx_k\|
= \tfrac{2}{1-\rho}\|x_k-\barx_k\|.
\end{align}
Therefore, from \cref{cor:avgin1}, \eqref{eq.error.bound.wk}, \eqref{eq.wk.dist.bound}, and~\eqref{eq:thm1-rho}, it follows that
\begin{align*}
    \|x_{k+1}-\barx\|^2&\leq \|x_k-\barx\|^2-\tfrac{1-\alpha}{\alpha}\|x_k-x_{k+1}\|^2\\
   & \le \|x_k-\barx\|^2-\tfrac{1-\alpha}{\alpha K_F^2}\|x_k-\barx_k\|^2\\
   & = \|x_k-\barx\|^2- (1-\rho^2)\|x_k-\barx_k\|^2\\
   & \le \|x_k-\barx\|^2-\tfrac{(1-\rho^2)(1-\rho)}{2}\|x_{k}-\barx\|^2.
\end{align*}
Equivalently,
\begin{align}
\label{eq:statement2}
    \|x_{k+1}-\barx\| &\le \left(1 - \tfrac{1+\rho}{2}(1-\rho)^2\right)^{\frac{1}{2}} \|x_k-\barx\|,
\end{align}
proving~\eqref{eq:thm1-hatrho}.

To finish, we prove the upper bounds for the rates $\rho$  and $\hat \rho$ in~\eqref{eq:thm1-hatrho}.
Let $\theta := \smash{\tfrac{1-\alpha}{\alpha K_F^2} }$. Clearly, $0 < \theta$ (as $0 < \alpha <1$), and by~\eqref{eq:positive} we obtain $\theta <1$. Now, we use
\begin{equation}\label{eq:implSimpl}
0 \le u < 1 \text{ implies } (1-u)^{\tfrac{1}{2}} \le 1 - \tfrac{1}{2}u.
\end{equation}
By~\eqref{eq:implSimpl} we obtain $\rho \le 1 - \tfrac{1}{2} \theta$, the bound for~$\rho$ in~\eqref{eq:thm1-rho}.
Since $f(x) := (1-x)^2$ is decreasing between $0 \le x < 1$, we have $f(\rho) \ge f(1-\tfrac{\theta}{2}) =  \tfrac {\theta^2}{4}$.
Using~\eqref{eq:implSimpl} we obtain, 
\[
\hat \rho = (1 - \tfrac{1+\rho}{2}f(\rho))^{\tfrac{1}{2}} \le 1 - \tfrac{1+\rho}{4}f(\rho) \le 1 -  \tfrac {\theta^2}{16},
\]
the bound for~$\hat \rho$ in~\eqref{eq:thm1-hatrho}.
\endproof

\cref{thm.FPlinconvergence} allows us to recast the problem of proving and estimating the rate of linear convergence as the problem of computing the constant $K_F$ in the generic error-bound condition~\eqref{eq.error.bound.F}.
In Section~\ref{sec:opt}, we will illustrate this by  estimating in forthright fashion,  the rate of linear convergence of different optimization algorithms to minimize well-known classes of functions.

By looking closely to the proof of \cref{thm.FPlinconvergence},  one can see that we have proven a stronger statement than the linear convergence of the FPI. We prove that if the error-bound condition~\eqref{eq.error.bound.F} holds for $R>0$ then the distance to the set of fixed-points decreases linearly when $F$ is applied to any point in $\{x: \dist(x,\Fix_F) \le R\}$.
Next, in \cref{prop.lcEB} we prove the converse of this statement, showing that the error-bound condition is necessary to obtain this form of linear convergence.

\begin{proposition}[Linear convergence implies error bound]\label{prop.lcEB}
Let $F:\R^n \to \R^n$ be such that $\Fix_F \ne \emptyset$. Let $R > 0$
and $\rho_F < 1$ be such that
\begin{equation}\label{eq.lc}
 \dist(x,\Fix_F) \le R \text{ implies }    \dist(F(x),\Fix_F) \le \rho_F \dist(x,\Fix_F).
\end{equation}
Then, the  error-bound condition~\eqref{eq.error.bound.F} holds for $R$, with constant $K_F = \smash{\tfrac 1{1-\rho_F}}$
\end{proposition}
\proof{}
Let $x$ be such that $\dist(x,\Fix_F) \le R$ holds. Let $\barx, \hat x \in \Fix_F$ be such that $\dist(x,\Fix_F) = \|x - \barx\|$ and $\dist(F(x),\Fix_F) = \|F(x) - \hat x\|$. Then,
\begin{align*}
\label{eq:num}
\|x - \bar x\|  & \le \|x - \hat x\| \\
& \le \|x-F(x)\| + \|F(x) - \hat{x}\| && \\
		  & \le \|x-F(x)\|  + \rho_F\|x- \bar{x}\| && (\text{by~\eqref{eq.lc}})
\end{align*}
Rearranging terms, we obtain,
\begin{equation}
\label{eq:error}
\dist(x,\{x:F(x) = x\}) = \|x - \bar x\| \le \frac{1}{1-\rho_F}\|x-F(x)\|.
\end{equation}
\endproof

\cref{thm.FPlinconvergence} and \cref{prop.lcEB} relate the rate of convergence $\rho_F$ of the FPI to the error-bound condition constant $K_F$. Using these results we can upper and lower bound on the convergence rate in terms of the error-bound constant and viceversa. We state these bounds in \cref{cor:ebc.lincon.bounds} and show that both sides of the bounds can be tight on Examples~\ref{ex:gradTight} and~\ref{ex:RotTight} below.

\begin{corollary}[Error-bound condition and linear convergence bounds]\label{cor:ebc.lincon.bounds}
Let $\alpha\in(0,1)$. Let $F:\R^n\to\R^n$ be a continuous $\alpha$-averaged operator with $\Fix_{F}\neq \emptyset$.
Let $R>0$ and define $\tilde{\rho}_F\in [0,1]$ and $\tilde{K}_F\in \R \cup \{\infty\}$ as
\[\tilde{\rho}_F = \sup_{x \notin \Fix_F}\frac{\dist(F(x),\Fix_F)}{\dist(x,\Fix_F) }
\text{ and }
\tilde K_F = \sup_{x \notin \Fix_F} \frac{\dist(x,\{x:F(x) = x\})} {\|F(x) - x\|},
\]
the tightest constants such that both the linear convergence result~\eqref{eq.lc}, and the error-bound condition~\eqref{eq.error.bound.F} hold.  We have the following equivalent relations,
\[1-\frac 1 {\tilde{K}_F}\le \tilde{\rho}_F \le \sqrt{1- \frac{1-\alpha}{\alpha \tilde{K}_F^2}}
\text{ and }
\sqrt{\frac{1-\alpha}{\alpha (1-\tilde{\rho}_F^2)}} \le \tilde{K}_F\le \frac 1{1-\tilde{\rho}_F}.\]
\end{corollary}

The next examples illustrate the tightness of the bounds on \cref{cor:ebc.lincon.bounds}.
\begin{example}[Convergence rates and error bounds for Gradient descent]\label{ex:gradTight}
Consider gradient descent applied to the optimization problem  $\min_x f(x) = \tfrac12 \|x\|^2$. Note that $\nabla f(x)=x$ is $1$-Lipschitz continuous. Let $\lambda\in(0,1)$. The gradient descent operator for this problem is given by $F_{GD}(x) = x - \lambda \nabla f(x)= (1-\lambda) x$.
In fact, $F_{GD}$ is a contraction with constant $\beta = 1-\lambda<1$ and unique fixed-point $\Fix_{F_{GD}} = \{0\}$. The fixed-point iteration of $F_{GD}$ converges at a global linear rate $\tilde\rho_{F_{GD}}=1-\lambda$ to this fixed-point. Next, we compute $\tilde K_{F_{GD}}$. For any $x\in\R^n$,
 \[
\dist(x,\Fix_{F_{GD}})=\|x-0\| =  \frac1\lambda \| x - F_{GD}(x)\|.
\]
Hence, $F_{GD}$ satisfies the global error bound with $\tilde K_{F_{GD}}=\frac1\lambda$. In particular, $\tilde \rho_{F_{GD}} = 1-\frac 1{\tilde K_{F_{GD}}}$.
\end{example}

\begin{example}[Convergence rates and error bounds for a rotation]\label{ex:RotTight}
Let $N_\theta:\R^2\to\R^2$ denote the rotation by the angle $\theta$ for some fixed $\theta \in (0,\pi)$. Let $F_\theta(x)= \tfrac12 (x+N_\theta(x))$ be the corresponding half-averaged operator. Notice that $0$ is the unique fixed point for $F_{\theta}$. Thus, for any $x \in \R^2$ we have $\dist(x,\Fix_{F_{\theta}}) = \|x\|$.

Standard trigonometric identities imply that, $\|F_\theta(x)\| = \cos(\tfrac12\theta)\|x\|$ and \( \|F_\theta(x)-x\| = \sin(\tfrac12\theta)\|x\|\).
Therefore,
\[
\tilde K_{F_\theta} = \sup_{x \in \R^2\setminus \{0\}} \frac{\dist(x,\Fix_{F_{\theta}})} {\|F_\theta(x) - x\|} = \sup_{x\in \R^2\setminus \{0\}} \frac{\|x\|} {\sin(\tfrac 12 \theta)\|x\|} = \frac1{\sin(\tfrac12\theta)}.
\]
On the other hand,
\[
\tilde \rho_{F_\theta} = \sup_{x \in \R^2\setminus \{0\}}\frac{\dist(F_\theta (x),\Fix_{F_\theta})}{\dist(x,\Fix_{F_\theta}) } = \sup_{x \in \R^2\setminus \{0\}}\frac{\|F_\theta(x)\|}{\|x\| } = \cos (\tfrac12 \theta).
\]
It follows that
\[
\tilde \rho_{F_\theta} = \cos (\tfrac12 \theta) = \sqrt{1- \sin^2(\tfrac12 \theta)} = \sqrt{1- \frac 1 {\tilde K_{F_\theta}^2}}.
\]
\end{example}

\section{Piecewise linear fixed-point operators.}
\label{sec:EB4PWLopers}
In order to specialize the results in \cref{thm.FPlinconvergence} and \cref{prop.lcEB} to
fixed-point operators associated with optimization algorithms, we now focus on
the class of piecewise linear fixed-point operators which, as discussed later in Section~\ref{sec:opt}, arise ubiquitously in optimization.

Specifically, we show in \cref{thm:HpieceWiseLin} that if an operator $F$ is piecewise linear, then $F$  satisfies the error-bound condition~\eqref{eq.error.bound.F} by showing that piecewise linear operator $I-F$ satisfies a {\em relative Hoffmann} bound~\citep{pena2021new}.
To show this formally, we use the following notation and terminology.

\begin{definition}[Compatible collection]
\label{def:compatible}
Let $P_1,\dots,P_k\subseteq \R^n$ be nonempty polyhedra and $F_i:\R^n\rightarrow \R^n, \;i=1,\dots,k$ be affine operators. Then we say that $(F_1,P_1),\dots,(F_k,P_k)$ is a {\em compatible collection} if the following two conditions hold:
\begin{itemize}\setlength{\itemindent}{1em}
    \item[(i)] $\R^n = P_1\cup\cdots\cup P_k$,
    \item[(ii)]  $F_i(x) = F_j(x)$ for all $i,j\in \{1,\dots,k\}$ and all $x \in P_i\cap P_j$.
\end{itemize}
\end{definition}
Given a compatible collection $(F_1,P_1),\dots,(F_k,P_k)$, we shall write $\bigcup_{i=1}^k F_i|P_i$ to denote the piecewise linear operator defined via
$x\in P_i \mapsto F_i(x)$.  Note that every piecewise linear operator $F$ is of the form $F=\bigcup_{i=1}^k F_i|P_i$ for some compatible collection $(F_1,P_1),\dots,(F_k,P_k)$.

In what follows, given an operator $G: \R^n \to \R^m$ we will denote the {\em set of  zeros} of the operator by $G^{-1}(0):= \{x \in \R^n: G(x) = 0\}$. \cref{thm:HpieceWiseLin} below relies on the following Hoffman-type constant.

\begin{definition}[Relative Hoffman constant~\citep{pena2021new}]
 Let $P\subseteq \R^n$ and $G:P\rightarrow \R^n$ be an affine operator on $P$ satisfying $G^{-1}(0)\cap P \ne \emptyset$.  Then Hoffman's Lemma~\cite{hoffman1952on,hoffman2003approximate} implies that
\begin{equation}\label{def.rel.Hoffman}
H(G|P):= \sup_{x\in P \setminus G^{-1}(0)} \frac{\dist(x,G^{-1}(0)\cap P)}{\|G(x)\|} < \infty
\end{equation}
with the convention that $H(G|P) = 0$ when $P\subseteq G^{-1}(0)$.
We call $H(G|P)$ the {\em Hoffman constant  of $G$ relative to $P$}. In the case when $P = \R^n$, we use the notation
\[
H(G) := H(G|\R^n);
\]
 that is $H(G)$ is the {\em standard} (i.e., non-relative) Hoffman constant.
\end{definition}

\begin{theorem}\label{thm:HpieceWiseLin}
Let $P_1,\dots,P_k\subseteq \R^n$ be nonempty polyhedra and $F_j:\R^n\rightarrow \R^n, \;j=1,\dots,k$ be affine operators such that $(F_1,P_1),\dots,(F_k,P_k)$ is a compatible collection. Let $F=\bigcup_{j=1}^k F_j|P_j$.
 If  $\Fix_{F} \neq \emptyset$ then the error-bound
 condition~\eqref{eq.error.bound.F} holds for the operator $F$ when $R$ small enough and
 \[
K_{F} = \max_{j \in \{1,\dots,k\}: P_j\cap {\Fix_F} \ne \emptyset} H(I-F_j|P_j).
\]
\end{theorem}
\proof{} Define the operators $G:= I - F$ and $G_j := I-F_j$ for $j=1,\dots,k$. Note that $\Fix_{F} \neq \emptyset$ implies that $G^{-1}(0) \neq \emptyset$. Let $H := \max_{j: P_j\cap G^{-1}(0)\neq \emptyset} H(G_j|P_j)$.  Proving the statement of the Lemma is equivalent to showing that there exists
$R>0$  such that for all $x \in \R^n$
\begin{equation}\label{eq:lemma2Equiv}
 \dist(x,G^{-1}(0))\le R \text{ implies }
\dist(x,G^{-1}(0)) \le H \|G(x)\|.
\end{equation}
To prove this, first consider  any $1 \le j \leq k$ such that $G_j^{-1}(0)\cap P_j \neq \emptyset$. Let $x \in P_j$. We have by~\eqref{def.rel.Hoffman},
\begin{equation}\label{eq.Hdist}
\dist(x,G^{-1}(0))\le \dist(x,P_j \cap G^{-1}(0)) \le H(G_j|P_j) \cdot \|G_j(x)\| \le H \cdot \|G(x)\|.
\end{equation}
Now, consider any $1 \le j \le k$ such that $G_j^{-1}(0)\cap P_j = \emptyset$. We have
$\epsilon_j := \inf_{x \in P_j} \|G_j(x)\| >0 $. Taking $R < H\epsilon_j$ we have,
\begin{equation}\label{eq.epsdist}
x \in P_j \text{ and } \dist(x,G^{-1}(0)) \le R \text{ implies }\dist(x,G^{-1}(0)) \le \tfrac R{\epsilon_j} \cdot \|G_j(x)\| \le H \cdot \|G(x)\|.
\end{equation}

Combining~\eqref{eq.Hdist} and~\eqref{eq.epsdist} we obtain equation~\eqref{eq:lemma2Equiv} for any $R \le H \cdot \min_{j: P_j\cap G^{-1}(0) \neq \emptyset}\epsilon_j$.
\endproof

\section{Piecewise linear-quadratic optimization.}
\label{sec:opt}

In this section, we leverage our main results, specifically \cref{thm.FPlinconvergence}  in combination with \cref{thm:HpieceWiseLin}, to demonstrate that a wide range of optimization algorithms, when applied to the class of \emph{piecewise linear-quadratic} (PLQ) optimization problems, satisfy the error-bound condition and achieve linear convergence (see \cref{prop:linconv}). The algorithms considered include gradient descent, proximal point, proximal gradient, Douglas-Rachford, and ADMM.

By PLQ optimization, we refer to the class of problems where the objective is a piecewise linear-quadratic convex function~\citep[see, e.g.,][Def. 10.20]{rockafellar2009variational} whose feasible set is a polyhedron. This class encompasses various fundamental types of optimization problems, including linear optimization (LO) and convex quadratic optimization (QO). Additionally, it includes least squares and its variations, such as LASSO and elastic net, as well as support vector machines (SVM).

Furthermore, we refine our results for the ADMM (equivalently Douglas-Rachford) algorithm to establish error bounds and linear convergence rates for its application to LO (\cref{sec:LP}) and QO  (\cref{sec:QP}) problems. For LO,  we obtain an error-bound constant \(K_{F_{DR}} =  1\), which translates into a linear rate of convergence \(\rho_{F_{DR}} \le \tfrac 12\), independent of the problem's data (\cref{thm:LP}). For QO, we derive an error-bound constant  and a linear rate of convergence, both of which depend solely on the condition number of the data matrix defining the problem's quadratic objective (\cref{thm:QP}).

To derive these results, we will first introduce some relevant notation and definitions.

\subsection{Preliminaries.}
\label{sec:prelim}
In what follows, we use standard notation and results from~\citep{boyd2004convex, rockafellar2009variational, rockafellar1997convex}.
We denote the set of  extended real numbers by $\eR := \R \cup \{+\infty\}$.
The identity matrix is denoted by $I$.
For any $A\in\R^{m\times n}$, we denote by $\rank(A)$ the column rank of $A$. We say $A$ has full rank if  $\rank(A)=n$.
For any $A\in\R^{m\times n}$ with $\rank(A\tr)=m$, we denote the Moore-Penrose or pseudo inverse as $A^\dag := A\tr(AA\tr)^{-1}$. For any proper differentiable function $f:\R^n \to \R$, the gradient of $f$, denoted $\nabla f$, is the vector of partial derivatives of $f$; that is $\nabla f(x) = (\tfrac{\partial f}{\partial x_1}(x), \cdots, \tfrac{\partial f}{\partial x_n}(x))$, and its Fenchel dual (or Fenchel conjugate), denoted $f^*$, is given by $f^*(x) = \sup_{u \in \R^n} ( x\tr u - f(u))$.

Let $\mu \ge 0$.
A function $f:\R^n \to \eR$ is $\mu$-strongly convex if $f(x)-\tfrac{\mu}{2}\|x\|^2$ is convex.
In particular, convex functions are 0-strongly convex.
We denote the class of proper, closed, $\mu$-strongly convex functions, $f:\R^n \to \eR $, by $\Gamma_\mu$.

The {\em proximal} operator is a fundamental concept in optimization that generalizes the notion of a projection.
Formally, given $f \in \Gamma_\mu$, the proximal operator $\prox_{f}:\R^n\rightarrow \R^n$, is defined by
\begin{equation}
\label{def:prox}
\prox_{f}(x):=\argmin_{u \in \R^n} \left\{f(u) + \frac{1}{2} \|u-x\|^2\right\}.
\end{equation}
As it is usual, given $\gamma >0$, we will consider computing the proximal operator on the {\em scaled} function $\gamma f$, denoted $\prox_{\gamma f}: \R^n\rightarrow \R^n$, which using~\eqref{def:prox} is equivalent to
\begin{equation}
\label{def:proxscale}
\prox_{\gamma f}(x)=\argmin_{u \in \R^n} \left\{\gamma f(u) + \frac{1}{2} \|u-x\|^2\right\}.
\end{equation}
If $f \in \Gamma_0$ then the proximal mapping is defined everywhere and the following {\em Moreau decomposition} holds. For any $x\in \R^n$
\begin{equation}\label{eq.moreau}
\prox_{\gamma f}(x) + \gamma\prox_{\tfrac{1}{\gamma} f^*}(\tfrac{1}{\gamma} x) = x.
\end{equation}
Let $\delta_S$ denote the {\em indicator function} 
over a closed and convex set $S\subseteq \R^n$. We have,
\[
\prox_{\delta_S}(x)=\Pi_S(x):= \argmin_{u \in S} \|u-x\|^2,
\]
the projection of $x$ onto $S$.
The proximal operator lays the foundation for proximal methods, a class of algorithms encompassing, for instance, the Douglas-Rachford algorithm, and it is a standard tool for solving nonsmooth, constrained, large-scale, or distributed optimization problems.

To simplify the presentation, it is  useful to introduce the \textit{reflector}  operator. Formally, given $\gamma >0$ and $f \in \Gamma_\mu$, the reflector operator $\refop_{\gamma f}: \R^n \to \R^n$ of the scaled function $\gamma f$ is defined by
\begin{equation}
\label{def:refop}
\refop_{\gamma f}(x) = 2\prox_{\gamma f}(x) - x,
\end{equation}
or equivalently by the equation
\begin{equation}
\label{def:refop2}
\prox_{\gamma f}(x) = \tfrac 12( x + \refop_{\gamma f}(x)).
\end{equation}

Both the proximal and the reflector operators are nonexpansive operators. To see this, it is useful to introduce the notion of {\em subdifferential}. The subdifferential of $f\in \Gamma_0$, at $x \in \R^n$, is the set  $\partial f \subset \R^n$ of all subgradients of $f$ at $x$. That is
\[
\partial f (x) = \{u \in \R^n: f(y)-f(x) \ge u\tr(y-x) \text{ for all }y \in \R^n\}.
\]
Let $\mu \ge 0$. For all $f \in \Gamma_\mu$, the subdifferential of $f$ is {\em monotone}. 
That is, for all
\begin{equation}\label{eq:diffMon}
(v-u)\tr(y-x) \ge 0 \text{ for all } x,y\in\R^n,\, u\in \partial f(x) \text{ and }v\in \partial f(y).
\end{equation}
The proximal operator~\eqref{def:proxscale}  is related to the subdifferential as follows
\begin{equation}\label{eq:prox->diff}
u = \prox_{\gamma f}(x)  \text{ if and only if }  \tfrac{1}{\gamma}(x-u)\in \partial f(u).
\end{equation}
This relationship allows us to characterize the fixed-points of the proximal operator. Namely, from~\eqref{eq:prox->diff}, it follows that
\begin{equation}
\label{eq:proxopt}
x \in \Fix_{\prox_{\gamma f}} \text{ if and only if } 0 \in \partial f(x)  \text{ if and only if } x \in \argmin f(x).
\end{equation}
Using~\eqref{def:refop2}, the condition~\eqref{eq:prox->diff} can be stated in terms of the reflector operator. That is
\begin{equation}\label{eq:refop->diff}
u = \refop_{\gamma f}(x)  \text{ if and only if }  \tfrac {1}{2 \gamma} (x-u) \in \partial f\left(\tfrac 12(x+u)\right).
\end{equation}
Let $x, y \in \R^n$ be given. From~\eqref{eq:refop->diff} and the monotonicity of the subdifferential~\eqref{eq:diffMon} we obtain
\(\tfrac{1}{2\gamma}(x-\refop_{\gamma f}(x))-\tfrac{1}{2\gamma}(y-\refop_{\gamma f}(y)))\tr (\tfrac{1}{2}(x+\refop_{\gamma f}(x))-\tfrac{1}{2}(y+\refop_{\gamma f}(y))) \ge 0
\) or, equivalently,
\begin{equation}\label{eq:refopnonexpansive}
\|\refop_{\gamma f}(x) - \refop_{\gamma f}(y)\|^2 \le  \|x - y\|^2.
\end{equation}
That is, the reflector is nonexpansive. Further, from~\eqref{def:refop}, it follows that $\prox_{\gamma f}$ is $\tfrac 12$-averaged.

\subsubsection{Douglas-Rachford Algorithm.}
\label{sec:DR}
Consider the optimization problem
\begin{equation}\label{eq.primal}
\dmin_{x \in \R^n} f(x) + g(x)
\end{equation}
where $f, g \in \Gamma_0$. Since the class $\Gamma_0$ includes indicator functions over convex sets, problem~\eqref{eq.primal} can be used to represent constrained convex optimization problems.
Algorithm~\ref{algo.DR} describes the Douglas-Rachford (DR) algorithm~\citep{douglas1956numerical,lions1979splitting, patrinos2014douglas} for solving problem~\eqref{eq.primal}.  

\begin{algorithm}[H]
{\small
  \caption{Douglas-Rachford} \label{algo.DR}
  \begin{algorithmic}[1]
    \State Pick $w_0\in \R^n$, $\gamma > 0$, $\alpha \in (0,1)$
\For{$k=0,1,2,\dots$} \label{step:two}
	\Statex 	\quad $x_{k+1}:=\prox_{\gamma f}(w_k)$
	\Statex 	\quad $y_{k+1}:=\prox_{\gamma g}(2x_{k+1}-w_{k})$
	\Statex     \quad $w_{k+1}:=w_k+2\alpha(y_{k+1}-x_{k+1})$
\EndFor
	\end{algorithmic}
	}
\end{algorithm}

Observe that the update in Step~\ref{step:two} of Algorithm~\ref{algo.DR} can also be written as the fixed-point iteration
\begin{equation}
\label{eq.update}
w_{k+1}:= F_{DR}(w_k),
\end{equation}
where the fixed-point {\em Douglas-Rachford} operator $F_{DR}:\R^n\rightarrow \R^n$ is defined as follows
\begin{equation}
\label{eq:FPI_DR}
F_{DR}(w) = w + 2\alpha (\prox_{\gamma g}(2\prox_{\gamma f}(w)-w) - \prox_{\gamma f}(w)).
\end{equation}
Then, using~\eqref{def:refop}, we have that
\begin{equation}
\label{eq:DRoperator}
F_{DR}(w) = (1 - \alpha) w  + \alpha\refop_{\gamma g}(\refop_{\gamma f}(w)).
\end{equation}
From~\eqref{eq:DRoperator}, we have that $F_{DR}$ is $\alpha$-averaged, as the reflector is nonexpansive and the class of nonexpansive operators is closed under composition. Therefore, by  \cref{thm:avgcon}, the DR algorithm converges sublinearly to a point in  $\Fix_{F_{DR}}$. Notice that $w \in \Fix_{F_{DR}}$ if and only $\refop_{\gamma g}(\refop_{\gamma f}(w)) = w$, which by~\eqref{eq:refop->diff} is equivalent to $\tfrac {1}{2 \gamma} (\refop_{\gamma f}(w)-w) \in \partial g(\tfrac 12 ( \refop_{\gamma f}(w) +w)) $. It also follows from~\eqref{eq:refop->diff}, that $\tfrac {1}{2\gamma} (w-\refop_{\gamma f}(w)) \in \partial f\left(\tfrac 12 (w + \refop_{\gamma f}(w))\right)$ and therefore,
\begin{equation}
\label{eq:DRtoopt}
w \in \Fix_{F_{DR}} \text{ is equivalent to }0 \in  \partial g\left(\tfrac 12 (w + \refop_{\gamma f}(w))\right)  + \partial f\left(\tfrac 12 (w + \refop_{\gamma f}(w))\right).
\end{equation}
As a result, the optimal solutions of~\eqref{eq.primal} can be characterized using the fixed-points of operator $F_{DR}$; by~\eqref{eq:DRtoopt} and~\eqref{def:refop2},
\begin{equation}
\label{eq:DRopt}
\argmin_{x \in \R^n} f(x) + g(x) = \{\prox_{\gamma f}(w): w \in \Fix_{F_{DR}}\}.
\end{equation}

\subsubsection{Peaceman-Rachford Algorithm.}
\label{sec:PR}
Note that if we purposely violate the domain restriction on $\alpha$ and set $\alpha = 1$ in Algorithm~\ref{algo.DR}, then the Douglas-Rachford algorithm becomes the Peaceman-Rachford algorithm~\citep[][Sec. 4]{combettes2011proximal}. Thus, the fixed-point operator associated to the Peaceman-Rachford algorithm is
\begin{equation}
\label{eq:PRoperator}
F_{PR}(w) = \refop_{\gamma g}(\refop_{\gamma f}(w)) = w + 2(\prox_{\gamma g}(2\prox_{\gamma f}(w)-w) - \prox_{\gamma f}(w)).
\end{equation}
From the discussion in \cref{sec:DR}, it follows that $F_{PR}$ is nonexpansive but not averaged. Further, it also follows that $w \in \Fix_{F_{PR}}$ if and only $\refop_{\gamma g}(\refop_{\gamma f}(w)) = w$. Thus, the optimality conditions~\eqref{eq:DRtoopt} and~\eqref{eq:DRopt} also hold for $F_{PR}$. However, our main results, \cref{thm.FPlinconvergence}, \cref{prop.lcEB}, and \cref{cor:ebc.lincon.bounds} can not be applied to the analysis of the Peaceman-Rachford algorithm.

\subsubsection{ADMM.}
\label{sec:ADMM}
Consider the optimization problem

\begin{equation}\label{eq.primaladmm}
\begin{array}{cl}
\dmin_{x \in \R^{n_1}, y \in \R^{n_2}}  & f(x) + g(y) \\
\st & Ax+By  = c,\\
\end{array}
\end{equation}
where $f:\R^{n_1} \to \R, g:\R^{n_2} \to \R \in \Gamma_0$, $A \in \R^{m \times n_1}$, $B \in \R^{m \times n_2}$, and $c \in \R^m$. Algorithm~\ref{algo.ADMM} describes the {\em scaled dual variable} ADMM algorithm~\citep[][Sec. 3.1.1]{boyd2011distributed}
for solving problem~\eqref{eq.primaladmm}.

\begin{algorithm}[H]
{\small
  \caption{ADMM for~\eqref{eq.primaladmm}} \label{algo.ADMM}
  \begin{algorithmic}[1]
    \State Pick $y_0 \in \dom(g)$, $w_0 \in \R^m$, $\rho > 0$
\For{$k=0,1,2,\dots$} \label{ADMMstep:two}
         \Statex \quad $x_{k+1} := \argmin_{x \in \R^{n_1}} \left (f(x) + w_k\tr (Ax + By_k-c) + \tfrac{\rho}{2}\|Ax+By_k-c\|^2 \right)$
	   \Statex \quad $y_{k+1}:=  \argmin_{y \in \R^{n_2}} \left (g(y) + w_k\tr (Ax_{k+1} + By-c) + \tfrac{\rho}{2}\|Ax_{k+1}+By-c\|^2 \right)$
	\Statex     \quad $w_{k+1}:= w_{k} + Ax_{k+1} +By_{k+1}-c$
\EndFor
	\end{algorithmic}
	}
\end{algorithm}
Now consider the ({\em negative}) Fenchel dual problem associated with~\eqref{eq.primaladmm}, that is
\begin{equation}
\label{eq:fencheldual}
\begin{array}{cl}
\dmin_{x \in \R^n}  &\hat{f}(x) +  \hat{g}(x),
\end{array}
\end{equation}
where $\hat{f}: \R^m \to \R$ and  $\hat{g}: \R^m \to \R$ are defined by
\begin{equation}
\label{eq:hats}
\hat{f}(x) = f^*(-A\tr x) + c\tr x,  \qquad
\hat{g}(x) = g^*(-B\tr x).
\end{equation}
It is well known that applying the ADMM algorithm to problem~\eqref{eq.primaladmm} is equivalent to applying the Douglas-Rachford algorithm, with $\alpha = \tfrac{1}{2}$ and $\gamma = 1/\rho$,  to  problem~\eqref{eq:fencheldual}~\citep[see, e.g.,][]{giselsson2016linear}. Using the definition of the DR operator~\eqref{eq:FPI_DR}, we obtain that the ADMM algorithm (\cref{algo.ADMM}) is equivalent to FPI of the  operator
\begin{equation}\label{def:FADDM}
F_{ADMM} (w):= w+\prox_{\gamma \hat{g}}(2\prox_{\gamma \hat{f}}(w)-w)-\prox_{\gamma \hat{f}}(w).
\end{equation}
\subsubsection*{Primal relationship between ADMM and DR.}
This {\em dual} relationship between ADMM and the DR algorithm leads to a {\em primal} relationship between the algorithms that we explain next. We will exploit this primal form in Sections~\ref{sec:algs},~\ref{sec:LP} and~\ref{sec:QP}.

Using the Moreau decomposition~\eqref{eq.moreau} and the reflector's definition~\eqref{def:refop}, it follows that for any $f,g \in \Gamma_0$, $\refop_{\gamma f}(x) = -\gamma \refop_{\tfrac{1}{\gamma} f^*}(\tfrac{1}{\gamma}x)$ and $\refop_{\gamma g}(x) = \gamma \refop_{\tfrac{1}{\gamma} g_*}(-\tfrac{1}{\gamma} x)$ where
$g_*(x) = g^*(-x)$. By~\eqref{def:FADDM}, we obtain then,
\begin{equation}
\label{eq:dualadmm}
F_{ADMM} (w)= \tfrac{1}{2}w  + \tfrac{1}{2}\refop_{\gamma  {\hat g}}(\refop_{\gamma {\hat f}}(w))=  \tfrac{1}{2} w  + \tfrac{\gamma}{2}\refop_{\tfrac{1}{\gamma}  {\hat g}_*}(\refop_{\tfrac{1}{\gamma}  {\hat f}^*}(\tfrac{1}{\gamma} w)).
\end{equation}
Computing the Fenchel duals of the functions $\hat f$ and $\hat g$, we obtain,
\begin{equation}
\label{eq:tildes}
\tilde{f}(w) :=  \hat f^*(w) = \min_{x \in \R^{n_1}} \{ f(x): c-Ax = w\} \text{ and }
\tilde{g}(w):=\hat g_*(w) = \min_{y \in \R^{n_2}} \{ g(y): By = w\}.
\end{equation}
It then follows from~\eqref{eq:dualadmm} and $\rho = \tfrac 1{\gamma}$, 
\begin{equation}
\label{eq:FADMMprimal}
\rho F_{ADMM}(\tfrac{1}{\rho} w) := \tfrac{1}{2}w  + \tfrac{1}{2}\refop_{\rho  {\tilde g} }(\refop_{\rho {\tilde f}}(w)),
\end{equation}
which is, the Douglas-Rachford operator corresponding to
\begin{equation}\label{eq:primalTildesDR}
\min_w \tilde f(w) + \tilde g(w)
\end{equation}
Notice that~\eqref{eq:primalTildesDR} is a equivalent reformulation of the primal problem~\eqref{eq.primaladmm}.

By~\citep[][Prop. 2.22]{rockafellar2009variational} $\tilde{f}, \tilde{g} \in \Gamma_0$. We obtain that  all the results in Section~\ref{sec:DR} apply to  $F_{ADMM}$. Also,
\begin{align*}
    \argmin_{(x, y) \in \mathbb{R}^{n_1} \times \mathbb{R}^{n_2}}
    &\left\{ f(x) + g(y) : Ax + By = c \right\} \\
= \; &\left\{
(\hat x,\hat y): \hat x=
\argmin_{x\in \R^{n_1}} \left\{ f(x) : c - Ax = \prox_{\gamma \tilde{f}}(w) \right\}, \right.  \\
&\qquad \left.
\hat y = \argmin_{y\in \R^{n_2}} \left\{ g(y) : By = \prox_{\gamma \tilde{f}}(w) \right\} , w \in \Fix_{F_{\mathrm{ADMM}}})
\right\}.
\end{align*}
The primal relation~\eqref{eq:FADMMprimal} between ADMM and DR implies also that Douglas-Rachford, when $\alpha = \tfrac 12$, can be seen as an application of ADMM to a particular setting.  Namely, consider the following reformulation of problem~\eqref{eq.primal} as a particular instance of problem~\eqref{eq.primaladmm},
\begin{equation}\label{eq.simpleadmm}
\begin{array}{cl}
\dmin_{x \in \R^n, y \in \R^n}  & f(x) + g(y) \\
\st & -x+y  = 0.\\
\end{array}
\end{equation}
In this particular case, $\tilde{f} = f$ and $\tilde{g} = g$. We obtain that
if $w_{k+1} = F_{ADMM}(w_k)$, then $\rho w_{k+1} = F_{DR}(\rho w_k)$.
This fact will be used in Section~\ref{sec:LP} and~\ref{sec:QP} to conclude that our results in those sections apply to both the ADMM algorithm and DR.

\subsubsection{From fixed-points to optimal points}
\label{sec:FPvsOP}

In this section we state the relation between the fixed-points of the operators in~\cref{tab:FPIs} and the corresponding optimal solutions of their associated optimization problems.
This relation was shown for the gradient descent fixed-point operator $F_{GD}$ in Example~\ref{ex:gradient}, for the proximal point fixed-point operator $F_{P}$ in~\eqref{eq:proxopt}, for the Douglas-Rachford $F_{DR}$ in \cref{sec:DR}, for the Peaceman-Rachford operator $F_{PR}$ in \cref{sec:PR} and for the ADMM operator $F_{ADMM}$ in \cref{sec:ADMM}.

\begin{lemma}[Relation between fixed points and optimality.]\label{lem:fix.p.opt.con}
Let $f,g \in\Gamma_0$, $\gamma > 0$. We have that,
\begin{enumerate}
    \item $\Fix_{F_{P}}=\{x\in \R^n: 0\in \partial f(x)\}=\argmin\{f(x):x \in \R^n\}$
    \item $\Fix_{F_{DR}}= \Fix_{F_{PR}}= \{y\in \R^n: 0\in \partial f(\prox_{\gamma f}(y))+ \partial g(\prox_{\gamma f}(y))\}$ and thus \\ $\argmin\{f(x)+g(x): x \in \R^n\}= \{\prox_{\gamma f}(y):y\in \Fix_{F_{DR}}\}$
    \item given $A \in \R^{m \times n_1}$, $B \in \R^{m \times n_2}$, and $c \in \R^m$, we obtain
    \[\Fix_{F_{\mathrm{ADMM}}} = \{w\in \R^m: 0\in \partial \tilde f(\prox_{\gamma\tilde f}(y))+ \partial \tilde g(\prox_{\gamma \tilde f}(y))\}.\]
     And,
    \begin{align*}
    \argmin_{(x, y) \in \mathbb{R}^{n_1} \times \mathbb{R}^{n_2}}
    &\left\{ f(x) + g(y) : Ax + By = c \right\} \\
= \; &\left\{
(\hat x,\hat y): \hat x=
\argmin_{x\in \R^{n_1}} \left\{ f(x) : c - Ax = \prox_{\gamma \tilde{f}}(w) \right\}, \right.  \\
&\quad \left.
\hat y = \argmin_{y\in \R^{n_2}} \left\{ g(y) : By = \prox_{\gamma \tilde{f}}(w) \right\} , w \in \Fix_{F_{\mathrm{ADMM}}}
\right\}.
\end{align*}
where  $\tilde{f}$, $\tilde{g}: \R^m \to \overline \R$ are defined by $\tilde{f}(w) :=  \min_{x \in \R^{n_1}} \{ f(x): c-Ax = w\}$ and $\tilde{g}(w):= \min_{y \in \R^{n_2}} \{ g(y): By = w\}$.
\end{enumerate}
If in addition, $f$ is differentiable, then
\begin{enumerate}[resume]
    \item $\Fix_{F_{GD}}=\{x\in\R^n: \nabla f(x)=0\}=\argmin\{f(x): x \in \R^n\}$
    \item $\Fix_{F_{PG}}=\{x\in \R^n: 0\in \nabla f(x)+ \partial g(x)\}=\argmin\{f(x)+g(x): x \in \R^n\}$
    \item given $S \subset \R^n$ a closed convex set,
     \[\Fix_{F_{GP}}=\{x\in S: \nabla f(x)\tr (y-x)\geq0 \text{ for all }y\in S\}=\argmin_{x\in S} f(x)\]
\end{enumerate}
\end{lemma}

\subsection{FPI algorithms for PLQ optimization.}
\label{sec:algs}

Under the assumptions that $f, g \in \Gamma_0$ are PLQ and $S$ is a polyhedron, all the optimization problems listed in Table~\ref{tab:FPIs} are PLQ optimization problems. Next, we argue that  as a result of \cref{thm.FPlinconvergence} and \cref{thm:HpieceWiseLin}, all the algorithms listed in Table~\ref{tab:FPIs} --- with the exception of the Peacemean-Rachford (PR) algorithm ---  have a linear rate of convergence to the optimal solution of their associated optimization problems and their associated fixed-point operator satisfies the error-bound condition~\eqref{eq.error.bound.F} for some $R >0$.

We start by providing a formal definition of PLQ functions.
\begin{definition}[{PLQ functions}]
A function $f: \R^n \to \eR$ is called {\em piecewise linear-quadratic} if there exist finitely many polyhedral sets $P_i \subseteq \R^n$, $i=1,\dots,k$, such that $\dom(f) = \bigcup_{i=1}^k P_i$ and on each $P_i$, $i=1,\dots,k$, $f$ is quadratic; that is,
for any $x \in P_i$, $f(x) = \tfrac{1}{2} x\tr Q_i x + c_i\tr x + b_i$ for some symmetric matrix $Q_i \in \R^{n \times n}$, $c_i \in \R^n$, and $b_i \in \R$, for $i=1,\dots,k$. If further, on each $P_i$, $i=1,\dots,k$, $f$ is affine, then $f$ is called {\em piecewise linear}.
\end{definition}

It \cref{prop:linconv} next, we show that for PLQ optimization the  algorithms in Table~\ref{tab:FPIs}, except for the PR algorithm, converge linearly to an optimal solution. Also, the corresponding fixed-point operators satisfy the error-bound condition. To obtain this result, first we argue that all these operators are averaged. Then, we show that when applied to PLQ optimization problems, these operators are piece-wise linear.
These two properties allow us to use \cref{thm.FPlinconvergence} and \cref{thm:HpieceWiseLin} to prove the desired convergence and error-bound
properties of the fixed-point operators of interest, which in turn correspond with the optimal solutions of their associated optimization problems (see~\cref{lem:fix.p.opt.con}).

\begin{proposition}[Linear convergence and error-bound condition for PLQ optimization]
\label{prop:linconv}
Assume that the functions $f,g \in \Gamma_0$ and PLQ, $\gamma > 0$, $\alpha \in (0,1)$, and $A, B \in \R^{n \times m}$, $c \in \R^m$.
\begin{enumerate}[label = (\roman*)]
\item Let $F \in \{F_P, F_{DR}, F_{ADMM}\}$, then $F$ satisfies the error-bound condition~\eqref{eq.error.bound.F} for some $R > 0$, and the sequence $x_k = F(x_{k-1})$ converges linearly to $x^* \in \Fix_F$, which corresponds to an
optimal solution of its associated optimization problem.
\item Let $F \in \{F_{GD}, F_{GP}, F_{PG}\}$. If additionally, $f$ is differentiable and either $\mu$-strongly convex for some $\mu >0$ or $L$-smooth for some $L>0$, and $S$ is a polyhedron, then $F$ satisfies the error-bound condition~\eqref{eq.error.bound.F} for some $R > 0$, and the sequence $x_k = F(x_{k-1})$ converges linearly to $x^* \in \Fix_F$, which corresponds to an
optimal solution of its associated optimization problem.
\end{enumerate}
\end{proposition}

\proof{} Let the fixed-point operator $F \in \{F_{GD}, F_{P}, F_{GP}, F_{PG}, F_{DR}, F_{ADMM}\}$ . Now we show that the assumptions imply that the operator $F$ is averaged.
Specifically, let $\gamma > 0$. As stated after~\eqref{eq:refopnonexpansive}, the fixed-point operator $F_{P}$ associated with the proximal algorithm is averaged. The fact that the fixed-point operators $F_{DR}$ and $F_{ADMM}$ respectively associated with the Douglas-Rachford and ADMM algorithms are averaged is shown in Sections~\ref{sec:DR} and~\ref{sec:ADMM} respectively. Further, under the additional assumption that $f$ is $L$-smooth for some $L > 0$ (resp. $\mu$-strongly convex for some $\mu > 0$), it was shown in Example~\ref{ex:gradient2} that  the fixed-point operator $F_{GD}$ associated with the gradient descent algorithm is $\tfrac{\lambda L}{2}$-averaged (resp. $\tfrac{\lambda \mu}{2}$-averaged) for any $0 <\lambda < \tfrac{2}{L}$
(resp. $0 <\lambda < \tfrac{2}{\mu}$). Using that the composition of averaged operators is averaged~\citep[][Prop. 2.4]{comb2014compositions}, it follows that  the  fixed-point operators $F_{GP}$ and $F_{PG}$, respectively associated with the gradient projection and the proximal gradient algorithms are also averaged.

Next, we show that when applied to PLQ optimization, all the fixed-point operators associated with the algorithms in Table~\ref{tab:FPIs} are piecewise linear. For this purpose, in the case of ADMM, it is relevant to note that since $f, g \in \Gamma_0$ and PLQ, then $\tilde{f}, \tilde{g}$ (defined in~\eqref{eq:tildes}) satisfy $\tilde{f}, \tilde{g} \in \Gamma_0$ and PLQ~\citep[][Prop. 11.32(b)]{rockafellar2009variational}.
Thus the piecewise linear property for all the operators follows from the fact that the gradient and proximal (and in particular the projection) operators map the class of  convex PLQ functions to the class of piecewise linear operators (\citep[][Prop. 10.21 and 12.30(d)]{rockafellar2009variational}).

Therefore, by \cref{thm:HpieceWiseLin}, it follows that there exists $R >0$ such that the error-bound condition~\eqref{eq.error.bound.F} holds. Thus, the linear convergence follows from \cref{thm.FPlinconvergence}, together with \cref{lem:fix.p.opt.con}.
\endproof

To conclude this section, we refine our results for the class of LO and QO problems. In particular we compute upper bounds on the error-bound constant for the Douglas-Rachford operator for these problem classes.
These bounds in turn provide bounds on the linear rate of convergence of the Douglas-Rachford algorithm for the same problem classes. Surprisingly our bounds do not depend on the conditioning of the constraint system of the feasible set.

\subsection{ADMM/DR for Linear Optimization.}\label{sec:LP}

Next, we analyze  the Douglas-Rachford algorithm and its corresponding fixed-point operator applied to LO. From the discussion in the last paragraph of \cref{sec:ADMM}, the ADMM algorithm and corresponding operator are equivalent to setting the parameter $\alpha = \tfrac 12$ in the Douglas-Rachford algorithm.

Let $m,n \in \mathbb{N}$ with $0 < n \le m$. Also, let $c \in \R^n$, $A \in \R^{m\times n}$, and $b \in \R^m$. Consider the problem
\begin{equation}\label{model.LP}
\begin{array}{rl}
\dmin_{x \in \R^n} & c\tr x\\
\text{s.t.}& Ax \le b,
\end{array}
\end{equation}
for which we assume an optimal solution exists; In particular, the feasible set $\{x \in \R^n: Ax \le b\} \neq \emptyset$. 

Problem~\eqref{model.LP} is a particular instance of the optimization problem~\eqref{eq.primal} with $g(x) = c\tr x$ and $f(x) = \delta_{X}(x)$, where \[X:= \{x\in\R^n:Ax\le b\}\] is the feasible region of the linear optimization problem~\eqref{model.LP}. Straightforward calculations  show that for any $x \in \R^n$,
  $\prox_{\gamma g}(x)=x-\gamma c$ and $\prox_{\gamma f}(x)=\Pi_X(x)$. Consequently, $\refop_{\gamma g}(x)=x - 2\gamma c$ and $\refop_{\gamma f}(x)=2\Pi_X(x)-x$. Thus, from~\eqref{eq:DRoperator}, it follows that fixed-point Douglas-Rachford operator for the linear optimization problem~\eqref{model.LP} is given by
\begin{equation}\label{eq.FDR.LP}
F_{DR}(x)= (1-2\alpha)x + 2\alpha \Pi_X(x)- 2\alpha \gamma c.
\end{equation}

For ease of exposition, in what follows, any $m \in \N$, we let $[m]:= \{1,\dots,m\}$
and for any $J \subseteq [m]$,  $A_J\in \R^{|J|\times n }$ denotes the submatrix of $A$ where the rows with indices in $J$ are selected, while $b_J\in \R^{|J|}$ is the vector constructed by selecting the corresponding entries from $b$.

We show that the error-bound condition \eqref{eq.error.bound.F} holds for $F_{DR}$ for some $R >0$.
To do this, we use that the projection onto a polyhedron is a piecewise linear operator. More precisely, given $J \subseteq [m]$, we let $X_J$ be the set of feasible points of~\eqref{model.LP} for which the linear constraints indexed by $J$ are tight; that is,
\[
X_J:= \{x \in X: A_J x = b_J\}.
\]
Also, let $P_J$ be the Minkowski sum of $X_J$ and the cone $A_J\tr \R^{|J|}_+$; that is,
\begin{equation}
\label{eq.PJ}
P_J:= X_J + A_J\tr \R^{|J|}_+.
\end{equation}
Further, notice that the projection of $u \in \R^n$ onto the polyhedron $\{x \in \R^n: A_Jx = b_J\}$, which we denote by $\Pi_J(u)$ for brevity, is given by
\[
\Pi_J(u):= (I - A_J^\dagger A_J) u + A^\dagger_J b_J.
\]
Notice that $P_{\emptyset}=X_{\emptyset} = X$ and for all $J \subseteq [m]$, if $X_J \neq \emptyset$, then $X_J$ is a face of $X$.
Moreover, $\R^n = \bigcup_{J \in \subseteq [m]: X_J \neq \emptyset} P_J$~\citep[][Ex. 6.16 and Thm. 6.46]{rockafellar2009variational}. Furthermore, if  $\rank(A_J\tr) < |J|$, there is $J' \subset J$ such that $\rank(A_{J'}\tr) = |J'|$ and $P_J = P_{J'}$.
Thus, defining $\cJ = \{J \subseteq [m]: \rank(A_J\tr) = |J|,\, X_J \neq \emptyset\}$ we have $\R^n=\bigcup_{J \in \cJ}P_J$.  
From Definition~\ref{def:compatible}, it follows that
that the nonempty polyhedra $P_{J}$ for all $J \in \cJ$ and the affine operators $\Pi_J$ for all $J \in \cJ$ form the compatible collection $\{(\Pi_J, P_J): J \in \cJ\}$.

 and for any $J \in \cJ$ and $x\in P_J$, $\Pi_X(x) = \Pi_{X_J}(x) = \Pi_J(x)$. Thus, we can write  the projection operator into the feasible set $X$ of~\eqref{model.LP} as $\Pi_X = \bigcup_{J \in \cJ} \Pi_{J}|P_J$. Therefore,
\begin{align*}
F_{DR}(x) = (1-2\alpha)x + 2\alpha \Pi_{J}(x) -2\alpha \gamma c = (I - 2\alpha A_J^\dagger A_J) x + 2 \alpha A^\dagger_J b_J - 2 \alpha \gamma c, \text{ for all }x \in P_J.
\end{align*}
Letting,
\begin{equation}\label{eq.GJ.LP}
G_J(x):=2\alpha (A_J^\dagger A_Jx-A_J^\dagger b_J +\gamma c ),
\end{equation}
we have
\begin{equation}
\label{eq:IminusFLP}
  I- F_{DR}=\bigcup_{J\in \cJ}G_J|P_J.
\end{equation}
That is, the piecewise linear operator $I-F_{DR}$ can be written in terms of the compatible collection of affine operators and nonempty polyhedra $\{(G_J, P_J):J \in \cJ\}$. Thus, specializing \cref{thm:HpieceWiseLin} we obtain an upper bound on the error-bound constant of $F_{DR}$, which using~\cref{thm.FPlinconvergence} translates into a bound on the convergence rate of the Douglas-Rachford algorithm. In particular, in \cref{thm:LP}, we obtain that the error-bound condition constant for the DR/ADMM operators  and the convergence rate of the DR/ADMM algorithms are independent of the problem data. Notice, however, that the value of~$R$ for which error-bound condition constant and the convergence rate are valid might depend on the problem data.

\begin{theorem}\label{thm:LP} Let $0<\alpha<1$. 
Assume \eqref{model.LP} has an optimal solution.  Then the error-bound condition~\eqref{eq.error.bound.F} holds for the Douglas-Rachford operator $F_{DR}$~\eqref{eq.FDR.LP} with $\sqrt{\tfrac{1-\alpha}{\alpha}} \le K_{F_{DR}} \le \tfrac{1}{2\alpha}$, when $R>0$ is sufficiently small. Consequently, the DR algorithm (Algorithm~\ref{algo.DR}) converges linearly with rate $\rho_{F_{DR}} \le 1-2\alpha(1-\alpha)$. In particular, for the ADMM operator $F_{ADMM}$~\eqref{eq:FADMMprimal} and the ADMM algorithm (Algorithm~\ref{algo.ADMM}) we obtain 
$K_{F_{ADMM}} = 1$ and $\rho_{F_{ADMM}} \le \tfrac 12$ respectively.
\end{theorem}

Before presenting the proof of \cref{thm:LP}, note that as a result of~\eqref{eq.GJ.LP} and~\eqref{eq:IminusFLP}, it follows from \cref{thm:HpieceWiseLin} that proving \cref{thm:LP} boils down to providing an appropriate bound on the Hoffman constant $H(G_J|P_J)$ for all $J \in \cJ$. Indeed, we will show that for all $J \in \cJ$, $H(G_J|P_J) \le H(G_J) \le \tfrac{1}{2\alpha}$. \cref{lem:subGzP} next shows that the first inequality in the previous equation holds.
\begin{lemma}\label{lem:subGzP} Let $J \subseteq \cJ$ be such that $G_J^{-1}(0)\cap P_J\neq \emptyset$, where $G_J(x)$ is defined by~\eqref{eq.GJ.LP} and $P_J$ by~\eqref{eq.PJ}.
Then \[H(G_J|P_J) \le H(G_J).\]
\end{lemma}
\proof{}
We claim that  for every $x\in P_J$,
\begin{equation}\label{eq:distLP}
    \dist (x, G_J^{-1}(0)\cap P_J) = \dist (x, G_J^{-1}(0)).
\end{equation}
By Definition~\eqref{def.rel.Hoffman},
\sloppy 
\[
H(G_J|P_J) = \sup_{x\in P_J\setminus G_J^{-1}(0)} \frac{\dist(x,G_J^{-1}(0)\cap P_J)}{\|G_J(x)\|}
= \sup_{x\in P_J\setminus G_J^{-1}(0)}  \frac{\dist(x,G_J^{-1}(0))}{\|G_J(x)\|} \le  H(G_J).
\]
To complete the proof we only need to prove~\eqref{eq:distLP}.
First, note that from~\eqref{eq.GJ.LP}, using $\Null(A_J^\dagger A_J) = \Null(A_J)$, we have that
\[
G_J^{-1}(0) := \{x \in \R^n:A_J^\dagger A_Jx = A_J^\dagger b_J - \gamma c\} = \{x \in \R^n:A_Jx = b_J - \gamma A_J c\}.
\]

For any $x \in \R^n$ the closest point in $G_J^{-1}(0)$ is then $\Pi_{G_J^{-1}(0)}(x) = (I - A_J^\dagger A_J) x + A_J^\dagger (b_J - \gamma A_J c)$. When $x \in P_J$, we have  $x = x'+ A_J\tr s$ with $x'\in X_J$ and $s \in \R_+^{|J|}$. Therefore
\begin{equation}
\label{eq:pj}
\Pi_{G_J^{-1}(0)}(x) = x' - A_J^\dagger (A_J x'- b_J) +    (I - A_J^\dagger A_J) A_J\tr s  - \gamma A_J^\dagger A_J  c = x'- \gamma A_J^\dagger A_J c.
\end{equation}
Moreover, from \eqref{eq.GJ.LP}, it also follows that $G_J^{-1}(0)\cap P_J \neq \emptyset$ implies
$\gamma c \in -A_J\tr \R_+^{|J|}$. Thus, for some $s' \in \R_+^{|J|}$, 
\[
\Pi_{G_J^{-1}(0)}(x) = x'+ A_J^\dagger A_J A_J\tr s' =   x'+  A_J\tr s' \in P_J.
\]
Therefore,
\[ \dist(x,G_J^{-1}(0) \cap P_J) = \dist(x,G_J^{-1}(0)).\]
\endproof

Now we prove \cref{thm:LP}.\par\medskip
\noindent 
{\em Proof of \cref{thm:LP}.}
 Let $J \subseteq [m]$ be such that $G_J^{-1}(0)\cap P_J\neq \emptyset$.
Let $M= 2\alpha A_J^\dagger A_J$ and $v = 2\alpha (A_J^\dagger b_J - \gamma c)$.
From~\eqref{eq.GJ.LP}, we have that   $G_J(x) = Mx-v$ for any $x\in \R^n$.
Then, by \cref{lem:subGzP},
\begin{align*}
H(G_J|P_J) & \le H(G_J) := \sup_{x\notin G_J^{-1}(0)} \frac{\dist(x,G_J^{-1}(0))}{\|G_J(x)\|}
 =  \sup_{x:Mx\neq v} \frac{\|M^\dagger(Mx - v)\|}{\|Mx - v\|} \\
 & \le \sup_{w\neq 0} \frac{\|M^\dagger w\|}{\|w\|} = \|M^\dagger\| =  \| (2\alpha A_J^\dagger A_J)^\dagger \| = \frac{1}{2\alpha} \| A_J^\dagger A_J \|.
\end{align*}
By considering the singular value decomposition $A_{J}$ one obtains that  $\|A_J^\dagger A_J\|=1$. Hence, we conclude that
$ H(G_J|P_J) \leq \frac{1}{2\alpha}$ for all $J \in \mathcal{J}$. Since~\eqref{eq.GJ.LP} and~\eqref{eq:IminusFLP} establish that the piecewise linear operator $I-F_{ADMM}$ for problem~\eqref{model.LP} can be expressed in terms of the compatible collection $\{(F_J, P_J):J \in \cJ\}$, it follows from \cref{thm:HpieceWiseLin} that the error-bound
 condition~\eqref{eq.error.bound.F} holds for the operator $F_{DR}$ when $R>0$ small enough and
$K_{F_{DR}} \le \frac{1}{2\alpha}$. The lower bound $K_{F_{DR}}$ follows from~\eqref{eq:thm1-rho} as $\rho_{F_{DR}} \ge 0$. This, together with the fact that $F_{DR}$ is averaged allow us to apply \cref{thm.FPlinconvergence} to obtain  $\rho_{F_{DR}} \le 1 - 2\alpha(1-\alpha)$. In particular the statement for ADMM follows by setting $\alpha = \tfrac{1}{2}$.
\endproof

\subsection{ADMM/DR for Quadratic Optimization.}\label{sec:QP}

Next, we analyze  the Douglas-Rachford algorithm and its corresponding fixed-point operator applied to QO. From the discussion in the last paragraph of \cref{sec:ADMM}, the ADMM algorithm and corresponding operator are equivalent to setting the parameter $\alpha = \tfrac 12$ in the Douglas-Rachford algorithm.

Throughout this section, we fix  $m,n \in \mathbb{N}$ with $0 < n \le m$. We also fix $Q\in\R^{n\times n}$ a non-zero symmetric positive semidefinite matrix, $c \in \R^n$, $A \in \R^{m\times n}$, and $b \in \R^m$.
Consider the linearly constrained (convex) quadratic optimization problem
\begin{equation}\label{model.QP}
\begin{array}{rl}
\dmin_x & \frac{1}{2}x\tr Q x + c\tr x\\
& Ax \le b,
\end{array}
\end{equation}
for which we assume a solution exists; that is $\{x \in \R^n: Ax \le b\} \neq \emptyset$. 

Problem~\eqref{model.QP} is a particular instance of~\eqref{eq.primal} with $g(x) = \frac{1}{2}x\tr Q x + c\tr x $ and $f(x) = \delta_{X}(x)$, where
\[
X:= \{x\in\R^n:Ax\le b\}
\]
is the feasible region of problem~\eqref{model.QP}. In this case,
  $\prox_{\gamma f}(x) = \Pi_X(x)$ and $\prox_{\gamma g}(x) = (\gamma Q+I)^{-1}(x-\gamma c)$. Consequently, $\refop_{\gamma f}(x) = 2\Pi_X(x)-x$ and $\refop_{\gamma g}(x) = 2(\gamma Q+I)^{-1}(x-\gamma c)-x$. Thus, from~\eqref{eq:DRoperator}, it follows that the fixed-point DR operator for the linearly constrained quadratic optimization problem~\eqref{model.QP} is given by
\begin{equation}\label{eq.FDR.QP}
F_{DR}(x)= x+2\alpha (\gamma Q+I)^{-1}\big(2\Pi_X(x)-x-\gamma c\big)-2\alpha \Pi_X(x).
\end{equation}

Thus, using the same tools and similar notation to the one introduced in Section~\ref{sec:LP}, 
it follows that for any $J \in \cJ$ and $x\in P_J$,
\begin{align}\label{eq.GJ.QP}
    I- F_{DR}(x)& =  G_J(x) :=  2 \alpha \Pi_{J}(x) - 2\alpha (\gamma Q+I)^{-1}\big(2\Pi_{J}(x)-x-\gamma c\big)\notag\\
    &= 2\alpha \big(I-A_J^\dagger A_J - (\gamma Q+I)^{-1}(I-2A_J^\dagger A_J)\big)x \notag\\ 
    & \qquad -2\alpha \big((\gamma Q+I)^{-1}(2A_J^\dagger b_J - \gamma c) - A_J^\dagger b_J\big)\notag\\
    &= M_J x - v_J,
\end{align}
where
\begin{equation}
\label{eq:Mdef}
    M_J = 2\alpha (I-A_J^\dagger A_J - (\gamma Q+I)^{-1}(I-2A_J^\dagger A_J)),
\end{equation}
and
\begin{equation}
\label{eq:vdef}
v_J = 2\alpha ((\gamma Q+I)^{-1}(2A_J^\dagger b_J - \gamma c) - A_J^\dagger b_J),
\end{equation}
and as in the linear programming case,
\begin{equation}
\label{eq:IminusFLPQP}
  I- F_{DR}=\bigcup_{J\in \cJ}G_J|P_J.
\end{equation}
That is, the piecewise linear operator $I-F_{DR}$ can be written in terms of the compatible collection of affine operators and nonempty polyhedra $\{(F_J, P_J):J \in \cJ\}$ defined above. Thus, specializing \cref{thm:HpieceWiseLin}, we can use \cref{thm.FPlinconvergence}. In particular, in \cref{thm:QP}, we obtain that the error-bound condition constant for the DR operator $F_{DR}$ and the convergence rate of the DR algorithm depends only on the condition number $\kQ := \tfrac{\lambda_{\max}(Q)}{\lambda_{\min}^+(Q)}$ of the matrix $Q$ defining the problem's quadratic objective (see \eqref{model.QP}), where $\lambda_{\min}^+(\cdot)$ is the smallest positive eigenvalue.
Notice, however, that the value of~$R$ for which error-bound condition and the convergence rate are valid might depend on the problem data.

\begin{theorem}\label{thm:QP}
Assume \eqref{model.QP} has a feasible solution. Let $0<\alpha<1$ and let $\gamma>0$ be such that  $\gamma_0 := \gamma \lambda_{\max}(Q)< 1$.
Then the error-bound condition~\eqref{eq.error.bound.F} holds for the Douglas-Rachford operator $F_{DR}$~\eqref{eq.FDR.QP} with
\begin{equation}
\label{eq:KDR_QP}
    K_{F_{DR}} \le \frac{1}{2\alpha} \frac{1+\gamma_0}{(1-\gamma_0)\gamma_0} (\kQ- \gamma_0),
\end{equation}
when $R>0$ is sufficiently small. Consequently, the Douglas-Rachford algorithm (Algorithm~\ref{algo.DR}) converges linearly with rate 
\begin{equation}
\label{eq:rhoDR_QP}
\rho_{F_{DR}} \le
1 - \frac{2\alpha(1 - \alpha)(1 - \gamma_0)^2 \gamma_0^2}{(1 + \gamma_0)^2  (\kQ - \gamma_0)^2}.
\end{equation}
In particular, for the ADMM operator $F_{ADMM}$~\eqref{eq:FADMMprimal} and the ADMM algorithm (Algorithm~\ref{algo.ADMM}) we obtain
\begin{equation*}
    K_{F_{ADMM}} \le \frac{1+\gamma_0}{(1-\gamma_0)\gamma_0} (\kQ- \gamma_0) \text{ and } \rho_{F_{ADMM}} \le
1 - \frac{(1 - \gamma_0)^2 \gamma_0^2}{2(1 + \gamma_0)^2  (\kQ - \gamma_0)^2}.
\end{equation*}
\end{theorem}

\begin{remark}[Compact expression for $K_{F_{DR}}$ and $\rho_{F_{DR}}$]
Set $\gamma$ such that $\gamma_0 = \tfrac 12$.  Then the upper bound on $K_{F_{DR}}$ in~\eqref{eq:KDR_QP} can be written as:
\[
  K_{F_{DR}} \le \tfrac{3}{\alpha}\kQ
\]
and then  the linear convergence rate $\rho_{F_{DR}}$ in~\eqref{eq:rhoDR_QP} satisfies
\[
 \rho_{F_{DR}} \le 1 - \frac{\alpha(1 - \alpha)}{18 \kQ^2}.
\]
\end{remark}

In order to prove \cref{thm:QP}, we first state an auxiliary result in \cref{lem:subGzP.QP}.

\begin{lemma}\label{lem:subGzP.QP}Let $J\in \cJ$. Let $M_J$ and $G_J$ be defined by~\eqref{eq.GJ.QP}-\eqref{eq:vdef}.
Assume $G_J^{-1}(0)\cap P_J\neq \emptyset$, where $P_J$ is given by~\eqref{eq.PJ}. Then 
\begin{enumerate}[label=(\roman*)]
\item $\Null(M_J) \subseteq \Null(Q) \cap \Null(A_J)$.
\label{it:subGzP.QP.one}
\item $G_J^{-1}(0)  \subseteq P_J$.
\label{it:subGzP.QP.two}
\end{enumerate}
\end{lemma}

\proof
Note that from the definition of $\cJ$, it follows that $A_J$ is full row rank; that is $\rank(A_J\tr) = |J|$. Let
$L_J := (\gamma Q+I)M_J = 2 \alpha (\gamma Q(I-A_J^\dagger A_J) + A_J^\dagger A_J)$,
and
$l_J := (\gamma Q+I)v_J = 2 \alpha((I-\gamma Q)A_J^\dagger b_J -\gamma c)$.

To prove~\ref{it:subGzP.QP.one}, let $w\in \Null(M_J) = \Null(L_J)$. Write $w = w^N+w^\perp$, for some $w^N\in \Null(A_J)$ and $w^\perp\in \Null(A_J)^\perp$. Then we have
\begin{equation}\label{eq:Hnull}
L_Jw^N = 2\alpha (\gamma Q(I-A_J^\dagger A_J) + A_J^\dagger A_J)w^N =  2\alpha \gamma Q w^N,
\end{equation}
and
\begin{equation}\label{eq:Hperp}
L_Jw^\perp = 2\alpha (\gamma Q(I-A_J^\dagger A_J) + A_J^\dagger A_J) w^\perp = 2\alpha w^\perp.
\end{equation}
The latter follows from the fact that $\Null(A_J)^\perp = \range(A_J\tr)$, and thus 
$A_J^\dagger A_J w^\perp = w^\perp$.

 Note that by \eqref{eq:Hnull} and \eqref{eq:Hperp} 
 \begin{equation}
 \label{eq:LJzero}
 0 = L_Jw = L_J(w^N+w^\perp)=2 \alpha(\gamma Qw^N+w^\perp). 
 \end{equation}
 Thus, $w^\perp = -\gamma Qw^N$. By the orthogonality of $w^N$ and $w^\perp$ we have
\[
0 = (w^N)\tr w^\perp  = - \gamma (w^N)\tr  Q w^N.
\]
By positive semidefiniteness of $Q$ we obtain $Qw^N=0$. By~\eqref{eq:LJzero}, $0 = L_Jw = 2 \alpha w^\perp$, and consequently $w^\perp=0$.
 Hence $w = w^N \in \Null(Q) \cap \Null(A_J)$. 
 
 To prove~\ref{it:subGzP.QP.two}, let $w_0 \in G_J^{-1}(0)\cap P_J$. As $\gamma Q+I$ is nonsingular, $G_J^{-1}(0) = \{ w \in \R^n: L_J w = l_J\} = w_0 + \Null(L_J)$. From~\ref{it:subGzP.QP.one}, $\Null(L_J) = \Null(M_J)  \subseteq \Null(A_J)$. Thus, $G_J^{-1}(0) \subseteq P_J + \Null(A_J) = P_J$.
\endproof

Now we provide the proof of \cref{thm:QP}, which is analogous to the proof of \cref{thm:LP}, the core of the proof is to bound the value of an appropriate Hoffman constant.\par\medskip
\noindent 
{\em Proof of \cref{thm:QP}.}
Let \( J \subseteq [m] \) be such that \( G_J^{-1}(0) \cap P_J \neq \emptyset \). Then 
\begin{align*}
H(G_J|P_J) &:= \sup_{w\in P_J\setminus G_J^{-1}(0)} \frac{\dist(w,G_J^{-1}(0)\cap P_J)}{\|G_J(w)\|} \\
 & \leq \sup_{w\notin G_J^{-1}(0)} \frac{\dist(w,G_J^{-1}(0))}{\|G_J(w)\|} & (\text{by \cref{lem:subGzP.QP}\ref{it:subGzP.QP.two}}) \notag \\
&= \frac 1{\sigma^+_{\min}(M_J)},
\end{align*}
where $\sigma^+_{\min}(\cdot)$ is the smallest positive singular value. 

To bound \( \sigma^+_{\min}(M_J) \), observe that:
\[
2\alpha I - (\gamma Q + I)M_J = 2\alpha (I - \gamma Q)(I - A_J^\dagger A_J).
\]
Since \( \gamma\lambda_{\max}(Q) < 1 \), the matrix \( I - \gamma Q \) is invertible, and we can write:
\[
I - A_J^\dagger A_J + \frac{1}{2\alpha}(I - \gamma Q)^{-1}(\gamma Q + I)M_J = (I - \gamma Q)^{-1}.
\]
Let $r=\rank(M_J)$, then $\sigma^+_{\min}(M_J) = \sigma_r(M_J)$. Using singular value bounds, we get:
\begin{equation}
\label{eq:bd1}
\sigma_{\max}(I - A_J^\dagger A_J) + \frac{1}{2\alpha} \sigma_{r}((I - \gamma Q)^{-1}(\gamma Q + I)M_J) \ge \sigma_{r}((I - \gamma Q)^{-1}) = \frac{1}{1 - \gamma \lambda_{r}(Q)}.
\end{equation}
From \cref{lem:subGzP.QP}\ref{it:subGzP.QP.one}, it follows that $r= \rank(M_J) \ge \rank(Q)$. Thus,
$\lambda_r(Q) \ge \lambda_{\min}^+(Q)$ which implies 
\begin{equation}
\label{eq:bd2}
\frac{1}{1 - \gamma \lambda_{r}(Q)} \le \frac{1}{1 - \gamma \lambda_{\min}^+(Q)} = \frac{\kQ}{\kQ - \gamma_0}.
\end{equation}
Since \( \sigma_{\max}(I - A_J^\dagger A_J) \le 1 \), it follows that:
\[
 \sigma_{r}((I - \gamma Q)^{-1}(\gamma Q + I)M_J) \ge \frac{2\alpha \gamma_0 }{\kQ - \gamma_0}.
\]
On the other hand, by sub-multiplicativity of singular values:
\begin{equation}
\label{eq:bd3}
\sigma_{r}((I - \gamma Q)^{-1}(\gamma Q + I)M_J)
\le 
\frac{1}{1 - \gamma \lambda_{\max}(Q)} \cdot (1 + \gamma \lambda_{\max}(Q)) \cdot \sigma_{r}(M_J).
\end{equation}
Combining~\eqref{eq:bd1}-\eqref{eq:bd3}
\[
\frac{2\alpha \gamma_0 }{\kQ - \gamma_0} \le  \frac{1 + \gamma_0}{1 - \gamma_0}  \sigma_{\min}^+(M_J).
\]
Therefore,
\[
H(G_J|P_J) \le \frac{1}{\sigma_{\min}^+(M_J)} \le \frac{1}{2\alpha} \frac{1+\gamma_0}{(1-\gamma_0)\gamma_0} (\kQ- \gamma_0).
\]

Since \( I - F_{DR} \) is affine over $P_J$, for each $J \in  \cJ$, \cref{thm:HpieceWiseLin} implies the error-bound condition holds for \( F_{DR} \) with:
\[
K_{F_{DR}} = \max_{J: P_J \cap \Fix_F \ne \emptyset} H(G_J|P_J) \le  \frac{1}{2\alpha} \frac{1+\gamma_0}{(1-\gamma_0)\gamma_0} (\kQ- \gamma_0).
\]
Finally, since \( F_{DR} \) is \( \alpha \)-averaged, \cref{thm.FPlinconvergence} implies~\eqref{eq:rhoDR_QP}.
In particular the statement for ADMM follows by setting $\alpha = \tfrac{1}{2}$.
\endproof

\section{Final Remarks.}
\label{sec:final}
Throughout the article, we have looked at some of the most popular optimization algorithms. However, it is clear that our results apply to an even wider range of optimization algorithms with similar properties when formulated in fixed-point iteration form. For example, consider the {\em plug-and-play} algorithms that are used in image processing~\citep{nair2021fixed,sinha2025linear}, in which new algorithms are created by replacing the proximal operator in the proximal gradient method or in the second update of ADMM by a simpler linear transformation referred as a {\em linear denoiser}. Under appropriate conditions on the linear denoiser, the fixed-point operator associated with these algorithms is averaged~\citep[see, e.g.,][Thm. 3.5 and 3.6]{nair2021fixed}. Thus our main results can be used to study the convergence of these algorithms.

Our results in Section~\ref{sec:QP} for quadratic optimization imply that one can scale the variables such that $Q$ is well-conditioned, without worrying about the effect on the conditioning of the linear constraints to improve the efficiency of the algorithm and the quality of the solution, as the convergence rates are independent of the data defining such linear constraints.

\printbibliography

\end{document}